\documentclass[oneside, reqno]{amsart}

\usepackage{xypic}
\input xy
\xyoption{all}
\usepackage{epsfig}
\usepackage{amsthm}
\usepackage{amssymb}
\usepackage{amsmath}
\usepackage{amscd}
\usepackage{color}
\usepackage[T1]{fontenc}
\usepackage{stackengine}
\usepackage{upgreek}
\usepackage[font=scriptsize]{caption}
\usepackage{wrapfig}
\stackMath

\newcommand{\bg}{\begin{equation}}
\newcommand{\ed}{\end{equation}}
\newcommand{\bga}{\begin{eqnarray}}
\newcommand{\eda}{\end{eqnarray}}
\newcommand{\pf}{\textbf{Proof:\ }}

\def\cbdu{\par{\raggedleft$\Box$\par}}

\newtheorem {Theorem}  {Theorem}

\numberwithin{Theorem}{section}

\newtheorem {Lemma}[Theorem]  {Lemma}

\theoremstyle{definition}
\newtheorem{Definition}[Theorem]{Definition}
\theoremstyle{remark}
\newtheorem{Remark}[Theorem]{\bf Remark}

\expandafter\chardef\csname pre amssym.def
at\endcsname=\the\catcode`\@ \catcode`\@=11
\def\undefine#1{\let#1\undefined}
\def\newsymbol#1#2#3#4#5{\let\next@\relax
 \ifnum#2=\@ne\let\next@\msafam@\else
 \ifnum#2=\tw@\let\next@\msbfam@\fi\fi
 \mathchardef#1="#3\next@#4#5}
\def\mathhexbox@#1#2#3{\relax
 \ifmmode\mathpalette{}{\m@th\mathchar"#1#2#3}%
 \else\leavevmode\hbox{$\m@th\mathchar"#1#2#3$}\fi}
\def\hexnumber@#1{\ifcase#1 0\or 1\or 2\or 3\or 4\or 5\or 6\or 7\or 8\or
 9\or A\or B\or C\or D\or E\or F\fi}

\font\teneufm=eufm10 \font\seveneufm=eufm7 \font\fiveeufm=eufm5
\newfam\eufmfam
\textfont\eufmfam=\teneufm \scriptfont\eufmfam=\seveneufm
\scriptscriptfont\eufmfam=\fiveeufm

\catcode`\@=\csname pre amssym.def at\endcsname

\newcounter{remark}
\def  \12  {{\frac{1}{2}}}

\def\build#1_#2^#3{\mathrel{\mathop{\kern 0pt#1}\limits_{#2}^{#3}}}

\numberwithin{equation}{section}

\begin{document}
%\currannalsline{0}{2006}

%\title[Local existence for Hall-MHD]{Local existence for the non-resistive Hall-magneto-hydrodynamics system in $\R^n$}
\title[Almost sure well-posedness for Hall MHD]{Almost sure well-posedness for Hall MHD}

%\author{hello}

\author [Mimi Dai]{Mimi Dai}

\address{Department of Mathematics, Statistics and Computer Science, University of Illinois at Chicago, Chicago, IL 60607, USA}
\address{School of Mathematics, Institute for Advanced Study, Princeton, NJ 08540, USA}
%\address{\color{white}...}
\email{mdai@uic.edu}

%\thanks{The author is grateful for the support of the NSF grants DMS--1815069 and DMS--2009422. }

\begin{abstract}
We consider the magnetohydrodynamics system with Hall effect accompanied with initial data in supercritical Sobolev space. 
%In particular, for initial data in  $\mathcal H^{s}(\mathbb T^n)\times \mathcal H^{s+1}(\mathbb T^n)$ with $s<\frac{n}{2}-1$, 
Via an appropriate randomization of the supercritical initial data,  
both local and small data global well-posedness for the system are obtained almost surely in critical Sobolev space.  
%$\mathcal H^{\frac{n}{2}-1}(\mathbb T^n)\times \mathcal H^{\frac{n}{2}}(\mathbb T^n)$.

%Goal: In particular, for initial data in  $\mathcal H^{s}(\mathbb T^n)\times \mathcal H^{s+1}(\mathbb T^n)$ with $s<\frac{n}{2}-1$, we randomize it in a suitable way.  We aim to show the existence and uniqueness of global in time solution for the system with such randomized data almost surely in space $\mathcal H^{\frac{n}{2}-1}(\mathbb T^n)\times \mathcal H^{\frac{n}{2}}(\mathbb T^n)$.

\bigskip


\vspace{0.05cm}
CLASSIFICATION CODE: 35Q35, 76D03, 76W05.
\end{abstract}

\maketitle

\section{Introduction}

A mathematical model for the incompressible magnetohydrodynamics (MHD) with Hall effect is given by
\begin{equation}\label{hmhd}
\begin{split}
u_t+(u\cdot\nabla) u-(B\cdot\nabla) B+\nabla \Pi=&\ \Delta u,\\
B_t+(u\cdot\nabla) B-(B\cdot\nabla) u+\nabla\times((\nabla\times B)\times B)=&\ \Delta B,\\
\nabla\cdot u=&\ 0
\end{split}
\end{equation}
on $\mathbb T^3\times [0,\infty)$. The unknowns $u$, $B$ and $\Pi$ are the velocity field, magnetic field and scalar pressure respectively. Note that $\nabla\cdot B(x,t)=0$ remains true for all the time $t>0$ if $\nabla\cdot B(x,0)=0$. The nonlinear term with the highest derivative is the Hall term $\nabla\times((\nabla\times B)\times B)$ which is posed to capture the rapid magnetic reconnection phenomena in plasma physics. The presence of the Hall term makes (\ref{hmhd}) a quasilinear system which is usually more challenging than semilinear systems. 
%With static background flow the Hall MHD system (\ref{hmhd}) reduces to the so-called electron MHD
%\begin{equation}\label{emhd0}
%\begin{split}
%B_t+\nabla\times((\nabla\times B)\times B)=&\ \Delta B,\\
%\nabla\cdot B=&\ 0.
%\end{split}
%\end{equation}

The well-posedness of the Hall MHD system (\ref{hmhd}) in various functional spaces has been studied extensively, for instance, see \cite{ADFL, CDL, CWW, Dai2, DT}. In these works the initial data and solution reside in the same spaces. In this paper, we are interested in the Cauchy problem for the Hall MHD with rough initial data yielding solutions with higher regularity. This can be achieved by randomizing the initial data properly. Such scheme for Cauchy problem with rough initial data was first developed in \cite{Bou96, BT1, BT2, LRS} for treating dispersive equations. It has been applied to the Navier-Stokes equation in \cite{NPS} to obtain global weak solutions when the initial data is in Sobolev space with negative index. Other applications can be found in \cite{CO, CD, Deng, DLM, Po} for both dispersive and dissipative systems. It is notable that the randomization strategy has the advantage to study Cauchy problems with supercritical feature, either the system being supercritical or the initial data being supercritical. 

In the author's previous work \cite{Dai2022}, we investigated the electron MHD with generalized diffusion $(-\Delta)^\alpha$ for suitable $\alpha>1$ in the supercritical regime. By randomizing the initial data in $H^s$ with $s<0$, we established global existence of weak solutions. When $\alpha>1$, the generalized electron MHD is no longer quasilinear, but semilinear. One major observation in \cite{Dai2022} is that there are obstructions to apply the approach of randomization of initial data for quasilinear equations, although the method is robust in analyzing supercritical equations. In this paper, we continue to study the full system of the Hall MHD with generalized diffusion,
\begin{equation}\label{hmhda}
\begin{split}
u_t+(u\cdot\nabla) u-(B\cdot\nabla) B+\nabla \Pi=&-(-\Delta)^\alpha u,\\
B_t+(u\cdot\nabla) B-(B\cdot\nabla) u+\nabla\times((\nabla\times B)\times B)=&-(-\Delta)^\alpha B,\\
\nabla\cdot u=&\ 0.
\end{split}
\end{equation}
We only consider $1\leq\alpha<\frac74$, in which regime system (\ref{hmhda}) is still supercritical, but semilinear rather than quasilinear. Unlike the classical MHD without Hall term, system (\ref{hmhda}) does not have a natural scaling. We extend the discussion on this aspect in the following.

When $B\equiv 0$, system (\ref{hmhda}) reduces to the hyperdissipative Navier-Stokes equation (NSE) 
\begin{equation}\label{hnse}
\begin{split}
u_t+(u\cdot\nabla) u+\nabla \Pi=&-(-\Delta)^\alpha u,\\
\nabla\cdot u=&\ 0
\end{split}
\end{equation}
which has the scaling property: if $(u(x,t), \Pi(x,t))$ solves (\ref{hnse}) with initial data $u_0(x)$, the rescaled pair $(u_\lambda, \Pi_\lambda)$ as
\begin{equation}\notag
u_\lambda(x,t)=\lambda^{2\alpha-1} u(\lambda x, \lambda^{2\alpha} t), \ \ \Pi_\lambda(x,t) = \lambda^{2(2\alpha-1)} \Pi(\lambda x, \lambda^{2\alpha}t)
\end{equation}
solves (\ref{hnse}) with initial data $u_{0,\lambda}=\lambda^{2\alpha-1}u_0(\lambda x)$. Among other scaling invariant (critical) spaces, the critical Sobolev space for (\ref{hnse}) in 3D is $\mathcal H^{\frac{5}{2}-2\alpha}$. Due to the prior energy estimate in $L^2=\mathcal H^0$, (\ref{hnse}) is critical for $\alpha=\frac54$; it is supercritical for $\alpha<\frac54$ and subcritical for $\alpha>\frac54$.

%System (\ref{emhd0}) has the natural scaling that if $B(x,t)$ solves the system with initial data $B_0(x)$, then the rescaled vector field
%\[B_\lambda(x,t)=B(\lambda x, \lambda^{2}t)\]
%for an arbitrary parameter $\lambda$ solves the system as well with initial data $B_0(\lambda x)$. Some scaling invariant spaces (also referred as critical spaces) with embedding for (\ref{emhd0}) in $n$-dimensional space are
%\begin{equation}\label{critical1}
%\dot{\mathcal H}^{\frac{n}2}\hookrightarrow L^{\infty}\hookrightarrow \dot B^{\frac{n}{p}}_{p,\infty} \hookrightarrow BMO \hookrightarrow \dot B^{0}_{\infty,\infty}, \ \ \ 1<p<\infty.
%\end{equation} 
%In view of (\ref{critical1}) we note the energy space $L^2(\mathbb T^n)$ is supercritical in both 2D and 3D, and system (\ref{emhd0}) is supercritical in both situations. Thus it is naturally challenging to analyze (\ref{emhd0}) even in 2D.

With static background flow $u\equiv 0$ the Hall MHD system (\ref{hmhd}) reduces to the so-called electron MHD
\begin{equation}\label{emhd}
\begin{split}
B_t+\nabla\times((\nabla\times B)\times B)=&-(-\Delta)^\alpha B,\\
\nabla\cdot B=&\ 0.
\end{split}
\end{equation}
System (\ref{emhd}) has the scaling
\begin{equation}\notag
B_\lambda(x,t)=\lambda^{2\alpha-2}B(\lambda x, \lambda^{2\alpha}t).
\end{equation}
The critical Sobolev space for (\ref{emhd}) on $\mathbb T^3$ is $\mathcal H^{\frac{7}2-2\alpha}$.
Again it follows from the basic energy law that (\ref{emhd}) is critical for $\alpha=\frac74$, supercritical for $\alpha<\frac74$ and subcritical for $\alpha>\frac74$.

In the full system (\ref{hmhda}), the scaling of the magnetic field equation plays a dominant role since it contains the highest degree nonlinear term. Thus it is supercritical when $\alpha<\frac74$. If $\alpha\geq \frac74$, global regular solution is expected for (\ref{hmhda}) through standard energy method. 
The purpose of the paper is to study the Cauchy problem of (\ref{hmhda}) in the regime $1\leq\alpha<\frac74$ and with initial data $(u_0, B_0)\in\mathcal H^{s_1}\times \mathcal H^{s_2}$, $s_1<\frac52-2\alpha$ and $s_2<\frac72-2\alpha$. By randomizing the initial data, we show the existence and uniqueness of solution to (\ref{hmhda}) in the space $\mathcal H^{\frac52-2\alpha}\times \mathcal H^{\frac72-2\alpha}$. As mentioned previously, the Hall term has the highest derivative in the system and poses the most challenges in our analysis. Therefore, we first investigate the hyperdissipative electron MHD (\ref{emhd}) and tackle the difficulties caused by the Hall term separately. We
then study the full coupled system (\ref{hmhda}), in which step the major difficulty comes from the coupling terms $(u\cdot\nabla) B$ and $(B\cdot\nabla) u$. The key ingredient to overcome the obstacles coming from the Hall term and coupling terms is the improved $L^p$ estimate for the free evolution of the randomized initial data. Details will be unfolded in later sections. We state the main results respectively for the electron MHD (\ref{emhd}), the hyperdissipative NSE (\ref{hnse}) and the full Hall MHD system (\ref{hmhda}) below.

\medskip

\begin{Theorem}\label{thm1}
Let $\alpha\in [1, \frac74)$. Let $B_0=f\in \mathcal H^{s}(\mathbb T^3)$ with $s\geq\max\{\frac{11}{2}-4\alpha, \frac{5}{2}-2\alpha\}$. There exists a set $\Sigma\subset \Omega$ with $P(\Sigma)=1$ such that for any $\omega\in\Sigma$ the electron MHD (\ref{emhd}) with the randomized initial data $f^\omega$ on $\mathbb T^3$ has a unique solution $B=B_{f^\omega}+H$ 
with $B_{f^\omega}=e^{-t(-\Delta)^\alpha} f^\omega$ and 
\begin{equation}\notag
H\in C([0,T]; \mathcal H^{\frac72-2\alpha}(\mathbb T^3))
\end{equation}
for some time $T>0$. Moreover, there exists a constant $c>0$ such that $T=\infty$ if $\|f\|_{\mathcal H^{s}(\mathbb T^3)}<c$.
\end{Theorem}

\medskip

\begin{Remark}
Note that $\frac{11}{2}-4\alpha\geq \frac{5}{2}-2\alpha$ for $1\leq\alpha\leq \frac32$, thus the theorem holds for
\begin{equation}\notag
\begin{split}
&s\geq \frac{11}{2}-4\alpha, \ \ \mbox{if} \ \ 1\leq\alpha\leq \frac32, \ \ \mbox{and}\\
&s\geq \frac{5}{2}-2\alpha, \ \ \mbox{if} \ \ \frac32<\alpha< \frac74.
\end{split}
\end{equation}
The gain of the derivative for the solution due to randomization of the initial data is
\begin{equation}\notag
\begin{cases}
(\frac72-2\alpha)-( \frac{11}{2}-4\alpha)=2\alpha-2, \ \ \mbox{if} \ \ 1<\alpha\leq \frac32, \\
(\frac72-2\alpha)-( \frac{5}{2}-2\alpha)=1, \ \ \mbox{if} \ \ \frac32<\alpha< \frac74.
\end{cases}
\end{equation}
In particular, when $\alpha\to 1$, the gain of the derivative is $2\alpha-2\to 0$. This corresponds to the classical well-posedness result in the critical space $\mathcal H^{\frac32}$ with initial data in the same space for the electron MHD with $\alpha=1$. Therefore the approach of randomization of the initial data seems not yield improvement in the quasilinear situation. 

\end{Remark}

%\medskip

\begin{Theorem}\label{thm-nse}
Let $\alpha\in [1, \frac74)$. Let $u_0=g\in \mathcal H^{s}(\mathbb T^3)$ with $s\geq\max\{\frac{7}{2}-\frac72\alpha, \frac{3}{2}-2\alpha\}$. There exists a set $\Sigma\subset \Omega$ with $P(\Sigma)=1$ such that for any $\omega\in\Sigma$ the NSE (\ref{hnse}) with the randomized initial data $g^\omega$ has a unique solution $u$ of the form 
\begin{equation}\notag
u=u_{g^\omega}+V
\end{equation}
%with  $u_{g^\omega}=e^{-t(-\Delta)^\alpha} g^\omega$, 
and 
\begin{equation}\notag
V\in C([0,T]; \mathcal H^{\frac52-2\alpha}(\mathbb T^3)). 
\end{equation}
for some time $T>0$. If in addition, $\|g\|_{\mathcal H^{s}(\mathbb T^3)}<c$ for some constant $c>0$, then $T=\infty$.
\end{Theorem}

%\medskip

\begin{Remark}
When $\alpha=1$, Theorem \ref{thm-nse} demonstrates almost sure well-posedness of the NSE in space $\mathcal H^{\frac12}(\mathbb T^3)$ for initial data in $L^2(\mathbb T^3)$. It recovers the result of \cite{CD}. 
\end{Remark}

%\medskip

\begin{Theorem}\label{thm2}
Let $\alpha\in [1, \frac74)$. Let $u_0=g\in \mathcal H^{s_1}(\mathbb T^3)$ and $B_0=f\in \mathcal H^{s_2}(\mathbb T^3)$ with $s_1\geq\max\{\frac{7}{2}-\frac72\alpha, \frac{3}{2}-2\alpha\}$ and $s_2\geq\max\{\frac{11}{2}-4\alpha, \frac{5}{2}-2\alpha\}$. There exists a set $\Sigma\subset \Omega$ with $P(\Sigma)=1$ such that for any $\omega\in\Sigma$ the Hall MHD (\ref{hmhda}) with the randomized initial data $(g^\omega, f^\omega)$ has a unique solution $(u,B)$ on $[0,T]$ for some $T>0$ in the form 
\begin{equation}\notag
u=u_{g^\omega}+V, \ \ B=B_{f^\omega}+H
\end{equation}
with  $u_{g^\omega}=e^{-t(-\Delta)^\alpha} g^\omega$ and $B_{f^\omega}=e^{-t(-\Delta)^\alpha} f^\omega$, and 
\begin{equation}\notag
V\in C([0,T]; \mathcal H^{\frac52-2\alpha}(\mathbb T^3)), \ \ H\in C([0,T]; \mathcal H^{\frac72-2\alpha}(\mathbb T^3)). 
\end{equation}
If in addition, $\|g\|_{\mathcal H^{s_1}(\mathbb T^3)}+\|f\|_{\mathcal H^{s_2}(\mathbb T^3)}<c$ for some constant $c>0$, then $T=\infty$.
\end{Theorem}

\begin{Remark}
A special case of Theorem \ref{thm2} is the almost sure well-posedness of (\ref{hmhda}) in the critical space $\mathcal H^{\frac52-2\alpha}\times \mathcal H^{\frac72-2\alpha}$ with initial data $(u_0, B_0)\in \mathcal H^{s-1}\times \mathcal H^{s}$ for $s\geq\max\{\frac{11}{2}-4\alpha, \frac{5}{2}-2\alpha\}$. When $\alpha=1$, that indicates almost sure well-posedness of the classical Hall MHD (\ref{hmhd}) in the critical space $\mathcal H^{\frac12}\times \mathcal H^{\frac32}$ with initial data $(u_0, B_0)\in \mathcal H^{\frac12}\times \mathcal H^{\frac32}$. In the deterministic case in \cite{Dai2}, well-posedness of (\ref{hmhd}) was established in nearly critical space $\mathcal H^{\frac12+\varepsilon}\times \mathcal H^{\frac32}$ with arbitrarily small $\varepsilon>0$. As pointed out in \cite{Dai2}, there is some essential obstacle to remove $\varepsilon$ due to the coupling feature of the Hall MHD system and different scalings for the velocity and magnetic field. The current result thus suggests that the randomization of initial data can be employed to remove such obstacle.   
\end{Remark}

\medskip

The rest of the paper is organized as: (i) in Section \ref{sec-pre} we lay out the notations to be used and recall the standard randomization procedure; (ii) Section \ref{sec-linear} is devoted to establishing improved estimates for the free evolution $u_{g^\omega}$ and $B_{f^\omega}$; (iii) we present a proof for Theorems \ref{thm1}, \ref{thm-nse} and \ref{thm2} in the last three sections respectively.

\bigskip

\section{Preliminaries}\label{sec-pre}

\medskip

\subsection{Notations}
We denote a general constant by $c$ which may differ in different estimates. Conventionally we denote $f\lesssim g$ for an inequality $f\leq cg$ with some constant $c>0$. For brevity, the Lebesgue space $L^p(\mathbb T^3)$ is sometimes denoted by $L^p_x$. It applies to Lebesgue spaces with respect to other variables as well,  for instance, $L^p(0,T)=L^p_t$. 

We introduce the weighted (in time) Lebesgue space
\begin{equation}\notag
L^{(r,s)}_tL^p_x=\left\{f: \left(\int_0^T \|t^r f(\cdot, t)\|_{L^p_x}^s \, dt\right)^{\frac1s}<\infty\right\}
\end{equation}
equipped with the norm 
\begin{equation}\notag
\|f\|_{L^{(r,s)}_tL^p_x}= \left(\int_0^T \|t^r f(\cdot, t)\|_{L^p_x}^s \, dt\right)^{\frac1s}.
\end{equation}

\medskip

\subsection{Randomization}

We recall the probabilistic estimates obtained in \cite{BT2}, which are valid for both Gaussian and Bernoulli variables. 

\begin{Lemma}\label{BT}
Let $(l_i(\omega))_{i=1}^\infty$ be a sequence of real-valued, zero-mean and independent random variables on a probability space $(\Omega, \mathcal A, P)$ with associated distributions $(\mu_i)_{i=1}^\infty$. Assume that there exists $c>0$ such that
\begin{equation}\notag%\label{ass1}
\left|\int_{-\infty}^{\infty}e^{\gamma x}\, d\mu_i(x)\right|\leq e^{c\gamma^2} \ \ \ \ \forall \gamma\in\mathbb R \ \ \ \forall i\geq 1.
\end{equation}
Then there exists $\beta>0$ such that 
\begin{equation}\notag
P\left(\omega: \left|\sum_{i=1}^\infty c_il_i(\omega)\right|>\lambda\right)\leq 2e^{-\frac{\beta \lambda^2}{\sum_{i=1}^\infty c_i^2}}   \ \ \ \ \forall \lambda>0 \ \ \ \forall (c_i)_{i=1}^\infty\in \ell^2.
\end{equation}
Consequently, there exists another constant $c>0$ such that 
\begin{equation}\notag
\left\|\sum_{i=1}^\infty c_il_i(\omega)\right\|_{L^q(\Omega)}\leq c \sqrt q\left(\sum_{i=1}^\infty c_i^2\right)^{\frac12}
 \ \ \ \ \forall q\geq 2 \ \ \ \forall (c_i)_{i=1}^\infty\in \ell^2.
\end{equation}
\end{Lemma}

\medskip

We adapt the standard diagonal randomization on the Sobolev space $\mathcal H^s(\mathbb T^n)$. 

\begin{Definition}\label{def-rand}
Let $(l_k(\omega))_{k\in \mathbb Z^n}$ be a sequence of real-valued and independent random variables on the probability space $(\Omega,\mathcal A, P)$ as in Lemma \ref{BT}. Let $e_k(x)=e^{ik\cdot x}$ for any $k\in \mathbb Z^n$. For a vector field $f=(f_1, f_2, ..., f_n)\in \mathcal H^s(\mathbb T^n)$ with Fourier coefficients $(a_{k})_{k\in \mathbb Z^n}$ and $a_k=(a_k^1, a_k^2, ..., a_k^n)$, the map 
\begin{equation}\label{R}
\begin{split}
\mathcal R:  (\Omega, \mathcal A) & \longrightarrow \mathcal H^s(\mathbb T^n)\\
\omega & \longrightarrow f^\omega, \ \ \ f^\omega(x)=\left(\sum_{k\in\mathbb Z^n}l_k(\omega)a_k^1 e_k(x), ... , \sum_{k\in\mathbb Z^n}l_k(\omega)a_k^n e_k(x)\right)
\end{split}
\end{equation}
equipped with the Borel sigma algebra is introduced. The map $\mathcal R$ is called randomization. 
\end{Definition}

It is worth to mention that the Leray projection $\mathbb P$ commutes with the randomization map $\mathcal R$, see \cite{NPS}. 
In view of Lemma \ref{BT}, we see that $\mathcal R$ is measurable and $f^\omega\in \mathcal H^s(\mathbb T^n)$ if $f\in  \mathcal H^s(\mathbb T^n)$. We also have
\begin{equation}\notag
f^\omega\in L^2(\Omega; \mathcal H^s(\mathbb T^n), \ \ \ \|f^\omega\|_{\mathcal H^s}\sim \|f\|_{\mathcal H^s}. 
\end{equation}
%Indeed, as shown in \cite{BT2}, the randomization $\mathcal R$ does not provide regularization of $\mathcal H^s$ in term of the regularity index $s$. Nevertheless, it gives rise to improved $L^p$ estimate almost surely. 

\bigskip

\section{Estimates of the free evolution with randomized initial data}\label{sec-linear}

Let $u(x,0)=g$ and $B(x,0)=f$ and their randomization $g^\omega$ and $f^\omega$ are given by
\begin{equation}\notag
\begin{split}
 g^\omega(x)=&\ \left(\sum_{k\in\mathbb Z^n}l_k(\omega)b_k^1 e_k(x), ... , \sum_{k\in\mathbb Z^n}l_k(\omega)b_k^n e_k(x)\right),\\
  f^\omega(x)=&\ \left(\sum_{k\in\mathbb Z^n}l_k(\omega)a_k^1 e_k(x), ... , \sum_{k\in\mathbb Z^n}l_k(\omega)a_k^n e_k(x)\right).
\end{split}
\end{equation}
We denote the free evolution of $g^\omega$ and $f^\omega$ according to the operator $-(-\Delta)^\alpha$ by $u_{g^\omega}=e^{-t(-\Delta)^\alpha}g^\omega$ and $B_{f^\omega}=e^{-t(-\Delta)^\alpha}f^\omega$ respectively. 

%We consider the linear equation with randomized initial data
%\begin{equation}\label{eq-lin}
%\begin{split}
%B_t+(-\Delta)^\alpha B=&\ 0,\\
%B(x,0)=&\ f^\omega,
%\end{split}
%\end{equation}
%and establish some estimates for its solution $B_{f^\omega}=e^{-t(-\Delta)^\alpha}f^\omega$. 

We first recall the Hardy-Littlewood-Sobolev lemma.

\begin{Lemma}\label{le-Hardy} Let $K(x)=|x|^{-\zeta}$ for $x\in\mathbb R^n$ and $\zeta>0$. Let $g\in L^p(\mathbb R^n)$. Then we have
\begin{equation}\notag
\|K* g\|_{L^q(\mathbb R^n)}\lesssim \|g\|_{L^p(\mathbb R^n)}
\end{equation}
with $1<p<q<\infty$ and $\zeta= n\left(1-\frac{1}{p}+\frac{1}{q}\right)$.
\end{Lemma}

The following two estimates will be used extensively in later sections. 

\begin{Lemma}\label{le-heat1}
Let $\alpha>0$, $p>0$ and $m>0$. Then
\begin{equation}\notag%\label{heat-lp}
\left\|e^{-t|\xi|^{2\alpha}}|\xi|^{m}\right\|_{L^p_{\xi}(\mathbb R^n)}\lesssim t^{-\frac{m}{2\alpha}-\frac{n}{2p\alpha}}.
\end{equation}
\end{Lemma}
\pf
Straightforward computation shows
\begin{equation}\notag
\begin{split}
\left\|e^{-t|\xi|^{2\alpha}}|\xi|^{m}\right\|_{L^p_{\xi}(\mathbb R^n)}=&\ \left(\int_{\mathbb R^n}e^{-pt|\xi|^{2\alpha}}|\xi|^{pm}\, d\xi\right)^{\frac1p}\\
=&\ t^{-\frac{m}{2\alpha}-\frac{n}{2p\alpha}}\left(\int_{\mathbb R^n}e^{-pt|\xi|^{2\alpha}}t^{\frac{pm}{2\alpha}}|\xi|^{pm}\, d(t^{\frac{1}{2\alpha}}\xi)\right)^{\frac1p}\\
\lesssim &\ t^{-\frac{m}{2\alpha}-\frac{n}{2p\alpha}}
\end{split}
\end{equation}
since the integral $\int_{\mathbb R^n} e^{-py^{2\alpha}} y^{pm}\, dy$ is bounded for $\alpha>0$, $p>0$ and $m>0$. It is obvious that the estimate holds on $\mathbb Z^n$ as well since $\mathbb Z^n\subset \mathbb R^n$.

\cbdu

\medskip

\begin{Lemma}\label{le-beta}
Let $0<r, s<1$ and $r+s=1$. There exists a constant $c>0$ independent of the time $t$ such that
\begin{equation}\notag
\int_0^t(t-\tau)^{-r}\tau^{-s}\, d\tau\leq c.
\end{equation}
\end{Lemma}
\pf
Changing variable $\tau= t\tau'$ in the integral gives
\begin{equation}\notag
\int_0^t(t-\tau)^{-r}\tau^{-s}\, d\tau=\int_0^1 (1-\tau')^{-r}(\tau')^{-s}\, d\tau' =B(1-s, 1-r)\leq c
\end{equation}
where $B(1-r, 1-s)$ is the Beta function.

\cbdu

\medskip
We establish some probabilistic estimates for $B_{f^\omega}$ and $u_{g^\omega}$ in the following.
\begin{Lemma}\label{le-heat2}
Let $r\geq q\geq p\geq 2$ and $s\geq -\frac{2\alpha}{q}+\frac{1}{2}-\frac{1}{p}$. The free evolution $B_{f^\omega}$ satisfies
\begin{equation}\notag
\|B_{f^\omega}\|_{L^r_\omega L^q_t L^p_x}\lesssim \sqrt r \|f\|_{\mathcal H^s_x}.
\end{equation}
%Consequently, if $-\frac{2\alpha}{q}+\frac{n}{2}-\frac{n}{p}\leq 0$, we have
%\begin{equation}\notag
%\|B_{f^\omega}\|_{L^r_\omega L^q_t L^p_x}\lesssim \sqrt r \|f\|_{L^2_x}.
%\end{equation}
\end{Lemma}
\pf
Recall 
\begin{equation}\notag
B_{f^\omega}(x,t)=e^{-t(-\Delta)^\alpha} f^\omega(x)=\sum_{k\in \mathbb Z^3} e^{-t|k|^{2\alpha}} l_{k}(\omega) a_{k} e_{k}(x).
\end{equation}
By Minkowski's inequality and Lemma \ref{BT}, we have
\begin{equation}\notag
\begin{split}
\|B_{f^\omega}\|_{L^r_\omega L^q_t L^p_x}\lesssim&\ \|B_{f^\omega}\|_{ L^q_t L^p_x L^r_\omega}\\
\lesssim&\ \sqrt r \left\|\left(\sum_{k\in \mathbb Z^n} e^{-2t|k|^{2\alpha}}  |a_{k}|^2 |e_{k}|^2 \right)^{\frac12}\right\|_{ L^q_t L^p_x}\\
=&\ \sqrt r \left\|\sum_{k\in \mathbb Z^n} e^{-2t|k|^{2\alpha}}  |a_{k}|^2 |e_{k}|^2 \right\|_{ L^{\frac{q}{2}}_t L^{\frac{p}{2}}_x}^{\frac12}.
\end{split}
\end{equation}
We continue the estimate with Minkowski's inequality again
\begin{equation}\notag
\begin{split}
\|B_{f^\omega}\|_{L^r_\omega L^q_t L^p_x}
\lesssim&\ \sqrt r \left\|\sum_{k\in \mathbb Z^n} e^{-2t|k|^{2\alpha}}  |a_{k}|^2 |e_{k}|^2 \right\|_{ L^{\frac{p}{2}}_x L^{\frac{q}{2}}_t}^{\frac12}\\
\lesssim&\ \sqrt r \left\|\sum_{k\in \mathbb Z^n} |k|^{-\frac{4\alpha}{q}}  |a_{k}|^2 |e_{k}|^2 \right\|_{ L^{\frac{p}{2}}_x}^{\frac12}\\
\lesssim&\ \sqrt r \left(\sum_{k\in \mathbb Z^n} |k|^{-\frac{4\alpha}{q}}  |a_{k}|^2 \|e_{k}\|^2_{ L^{p}_x}\right)^{\frac12}\\
\end{split}
\end{equation}
where we used the estimate
\begin{equation}\notag
\begin{split}
\|e^{-2t|k|^{2\alpha}}\|_{L^{\frac{q}{2}}_t}=&\ \left(\int_0^t e^{-qt|k|^{2\alpha}}\, dt\right)^{\frac{2}{q}}\\
=&\ |k|^{-\frac{4\alpha}{q}}\left(\int_0^t e^{-qt|k|^{2\alpha}}\, d(t|k|^{2\alpha})\right)^{\frac{2}{q}}\\
\lesssim&\ |k|^{-\frac{4\alpha}{q}}.
\end{split}
\end{equation}
Note that for the basis $\{e_k\}$ we have
% see the paper "Concerning the $L^p$ norm of spectral clusters for second order elliptic operators on compact manifolds", C. Sogge, J. Funct. Anal., 77: 123-138, 1988
\[\|e_{k}\|^2_{ L^{p}_x}\lesssim |k|^{2(\frac{1}{2}-\frac{1}{p})}\|e_{k}\|^2_{ L^{2}_x}.\]
Thus we further deduce 
\begin{equation}\notag
\begin{split}
\|B_{f^\omega}\|_{L^r_\omega L^q_t L^p_x}
\lesssim&\ \sqrt r \left(\sum_{k\in \mathbb Z^n} |k|^{-\frac{4\alpha}{q}}  |a_{k}|^2 \|e_{k}\|^2_{ L^{p}_x}\right)^{\frac12}\\
\lesssim&\ \sqrt r \left(\sum_{k\in \mathbb Z^n} |k|^{-\frac{4\alpha}{q}+1-\frac{2}{p}}  |a_{k}|^2 \|e_{k}\|^2_{ L^{2}_x}\right)^{\frac12}\\
\lesssim&\ \sqrt r \|f\|_{\mathcal H_x^{-\frac{2\alpha}{q}+\frac{1}{2}-\frac{1}{p}}}.
\end{split}
\end{equation}

\cbdu

\medskip

\begin{Lemma}\label{le-heat3}
Let $r\geq \frac{1}{2\beta}\geq p\geq 2$ and $\frac{1}{\alpha}+\frac{3}{2p\alpha}+2\beta=1$. Assume $s> \frac52-2\alpha$. 
%{\color{red}Double check in later computation to make sure it is $\frac{1}{\alpha}+\frac{3}{2p\alpha}+2\beta=1$ or $\frac{1}{\alpha}+\frac{1}{2p\alpha}+2\beta=1$.} 
The estimate
\begin{equation}\notag
\|B_{f^\omega}\|_{L^r_\omega L^{\frac{1}{2\beta}}_t L^p_x}\lesssim \sqrt r \|f\|_{\mathcal H^s_x}
\end{equation}
is valid.
\end{Lemma}
\pf
It is a special case of Lemma \ref{le-heat2}. Indeed, taking $q=\frac{1}{2\beta}$ in Lemma \ref{le-heat2}, we obtain 
\begin{equation}\notag
s\geq -\frac{2\alpha}{q}+\frac{1}{2}-\frac{1}{p}=-4\alpha\beta-\frac{1}{p}+\frac12=-2\alpha\left(2\beta+\frac{1}{2p\alpha}\right)+\frac12.
\end{equation}
On the other hand, we observe
\[2\beta+\frac{1}{2p\alpha}=1-\frac{1}{\alpha}-\frac{1}{p\alpha}<1-\frac{1}{\alpha}.\]
Therefore it follows that
\begin{equation}\notag
s\geq -2\alpha\left(2\beta+\frac{1}{2p\alpha}\right)+\frac12>-2\alpha \left(1-\frac{1}{\alpha}\right)+\frac12=\frac52-2\alpha.
\end{equation}
\cbdu

\medskip

\begin{Lemma}\label{le-heat4}
Let $r\geq q\geq p\geq 2$ and $\eta>0$. Then the estimate
\begin{equation}\notag
\|B_{f^\omega}\|_{L^r_\omega L^{(\eta, q)}_t L^p_x}\lesssim \sqrt r \|f\|_{\mathcal H^s_x}
\end{equation}
holds for $s\geq -2\alpha\eta-\frac{2\alpha}{q}+\frac{1}{2}-\frac{1}{p}$. 
\end{Lemma}
\pf
Thanks to Minkowski's inequality and Lemma \ref{BT} again, we infer
\begin{equation}\notag
\begin{split}
\|B_{f^\omega}\|_{L^r_\omega L^{(\eta, q)}_t L^p_x}
=&\ \left\| \sum_{k\in \mathbb Z^n} t^\eta e^{-t|k|^{2\alpha}} l_{k}(\omega) a_{k} e_{k}(x) \right\|_{L^r_\omega L^{q}_t L^p_x} \\
\lesssim&\ \sqrt r \left\|\left(\sum_{k\in \mathbb Z^n} t^{2\eta}e^{-2t|k|^{2\alpha}}  |a_{k}|^2 |e_{k}|^2 \right)^{\frac12}\right\|_{ L^q_t L^p_x}\\
=&\ \sqrt r \left\|\sum_{k\in \mathbb Z^n} t^{2\eta} e^{-2t|k|^{2\alpha}}  |a_{k}|^2 |e_{k}|^2 \right\|_{ L^{\frac{q}{2}}_t L^{\frac{p}{2}}_x}^{\frac12}\\
\lesssim&\ \sqrt r \left\|\sum_{k\in \mathbb Z^n} t^{2\eta} e^{-2t|k|^{2\alpha}}  |a_{k}|^2 |e_{k}|^2 \right\|_{ L^{\frac{p}{2}}_x L^{\frac{q}{2}}_t}^{\frac12}.
\end{split}
\end{equation}
The norm in time can be estimated as
\begin{equation}\notag
\begin{split}
\|t^{2\eta}e^{-2t|k|^{2\alpha}}\|_{L^{\frac{q}{2}}_t}=&\ \left(\int_0^t t^{q\eta}e^{-qt|k|^{2\alpha}}\, dt\right)^{\frac{2}{q}}\\
=&\ |k|^{-4\alpha\eta-\frac{4\alpha}{q}}\left(\int_0^t (t|k|^{2\alpha})^{q\eta}e^{-qt|k|^{2\alpha}}\, d(t|k|^{2\alpha})\right)^{\frac{2}{q}}\\
\lesssim&\ |k|^{-4\alpha\eta-\frac{4\alpha}{q}}.
\end{split}
\end{equation}
Therefore we continue the estimate as
\begin{equation}\notag
\begin{split}
\|B_{f^\omega}\|_{L^r_\omega L^{(\eta, q)}_t L^p_x}
\lesssim&\ \sqrt r \left\|\sum_{k\in \mathbb Z^n} t^{2\eta} e^{-2t|k|^{2\alpha}}  |a_{k}|^2 |e_{k}|^2 \right\|_{ L^{\frac{p}{2}}_x L^{\frac{q}{2}}_t}^{\frac12}\\
\lesssim&\ \sqrt r \left\|\sum_{k\in \mathbb Z^n} |k|^{-4\alpha\eta-\frac{4\alpha}{q}}  |a_{k}|^2 |e_{k}|^2 \right\|_{ L^{\frac{p}{2}}_x}^{\frac12}\\
\lesssim&\ \sqrt r \left(\sum_{k\in \mathbb Z^n} |k|^{-4\alpha\eta-\frac{4\alpha}{q}}  |a_{k}|^2 \|e_{k}\|^2_{ L^{p}_x}\right)^{\frac12}\\
\lesssim&\ \sqrt r \left(\sum_{k\in \mathbb Z^n} |k|^{-4\alpha\eta-\frac{4\alpha}{q}+1-\frac{2}{p}}  |a_{k}|^2 \|e_{k}\|^2_{ L^{2}_x}\right)^{\frac12}\\
\lesssim&\ \sqrt r \|f\|_{\mathcal H_x^{-2\alpha\eta-\frac{2\alpha}{q}+\frac{1}{2}-\frac{1}{p}}}.
\end{split}
\end{equation}

\cbdu

\medskip

\begin{Lemma}\label{le-heat5}
Let $r\geq \frac{1}{\beta}\geq p\geq 2$ and $\frac{1}{\alpha}+\frac{3}{2p\alpha}+2\beta=1$. Let $s> \frac52-2\alpha$. We have the estimate
\begin{equation}\notag
\|B_{f^\omega}\|_{L^r_\omega L^{(\beta, \frac{1}{\beta})}_t L^p_x}\lesssim \sqrt r \|f\|_{\mathcal H^s_x}.
\end{equation}
\end{Lemma}
\pf
This is a special case of Lemma \ref{le-heat4} with $\eta=\beta$ and $q=\frac{1}{\beta}$.
\cbdu

\medskip

\begin{Lemma}\label{le-heat6}
Let $r\geq 2$ and $0<\eta\leq\frac12$. Assume $s\geq \frac{11}{2}-4\alpha$. The estimate
\begin{equation}\notag
\|B_{f^\omega}\|_{L^r_\omega L^{(\eta, 2)}_t \mathcal H^{\frac{11}{2}-2\alpha}_x}\lesssim \sqrt r \|f\|_{\mathcal H^s_x}
\end{equation}
holds. For $s\geq \frac{11}{2}-\frac{7\alpha}2$, we have
\begin{equation}\notag
\|B_{f^\omega}\|_{L^r_\omega L^{(\eta, 4)}_t \mathcal H^{\frac{11}{2}-2\alpha}_x}\lesssim \sqrt r \|f\|_{\mathcal H^s_x}.
\end{equation}
\end{Lemma}
\pf
We only show details for the proof of the first inequality; the proof of the second inequality is analogous. 
It follows from Minkowski's inequality and Lemma \ref{BT} that
\begin{equation}\notag
\begin{split}
\|B_{f^\omega}\|_{L^r_\omega L^{(\eta, 2)}_t \mathcal H^{\frac{11}{2}-2\alpha}_x}
=&\ \|t^{\eta}(1-\Delta)^{\frac{11}{4}-\alpha}B_{f^\omega}\|_{L^r_\omega L^{2}_t L^2_x}\\
=&\ \left\| \sum_{k\in \mathbb Z^n} t^\eta (1+|k|^2)^{\frac{11}{4}-\alpha}e^{-t|k|^{2\alpha}} l_{k}(\omega) a_{k} e_{k}(x) \right\|_{L^r_\omega L^{2}_t L^2_x} \\
\lesssim&\ \sqrt r \left\|\left(\sum_{k\in \mathbb Z^n} t^{2\eta}(1+|k|^2)^{\frac{11}{2}-2\alpha}e^{-2t|k|^{2\alpha}}  |a_{k}|^2 |e_{k}|^2 \right)^{\frac12}\right\|_{ L^2_t L^2_x}\\
=&\ \sqrt r \left\|\sum_{k\in \mathbb Z^n} t^{2\eta} (1+|k|^2)^{\frac{11}{2}-2\alpha}e^{-2t|k|^{2\alpha}}  |a_{k}|^2 |e_{k}|^2 \right\|_{ L^{1}_t L^{1}_x}^{\frac12}\\
\lesssim&\ \sqrt r \left\|\sum_{k\in \mathbb Z^n} t^{2\eta} (1+|k|^2)^{\frac{11}{2}-2\alpha}e^{-2t|k|^{2\alpha}}  |a_{k}|^2 |e_{k}|^2 \right\|_{ L^{1}_x L^{1}_t}^{\frac12}.
\end{split}
\end{equation}
We estimate the norm in time as before
\begin{equation}\notag
\begin{split}
\|t^{2\eta}e^{-2t|k|^{2\alpha}}\|_{L^{1}_t}=&\ \int_0^t t^{2\eta}e^{-qt|k|^{2\alpha}}\, dt\\
=&\ |k|^{-4\alpha\eta-2\alpha}\int_0^t (t|k|^{2\alpha})^{2\eta}e^{-qt|k|^{2\alpha}}\, d(t|k|^{2\alpha})\\
\lesssim&\ |k|^{-4\alpha\eta-2\alpha}.
\end{split}
\end{equation}
Thus we infer
\begin{equation}\notag
\begin{split}
\|B_{f^\omega}\|_{L^r_\omega L^{(\eta, 2)}_t \mathcal H^{\frac{11}{2}-2\alpha}_x}
\lesssim&\ \sqrt r \left\|\sum_{k\in \mathbb Z^n} t^{2\eta} (1+|k|^2)^{\frac{11}{2}-2\alpha}e^{-2t|k|^{2\alpha}}  |a_{k}|^2 |e_{k}|^2 \right\|_{ L^{1}_x L^{1}_t}^{\frac12}\\
\lesssim&\ \sqrt r \left\|\sum_{k\in \mathbb Z^n} |k|^{-4\alpha\eta-2\alpha+11-4\alpha}  |a_{k}|^2 |e_{k}|^2 \right\|_{ L^{1}_x}^{\frac12}\\
\lesssim&\ \sqrt r \left(\sum_{k\in \mathbb Z^n} |k|^{-4\alpha\eta-6\alpha+11}  |a_{k}|^2 \|e_{k}\|^2_{ L^{2}_x}\right)^{\frac12}\\
\lesssim&\ \sqrt r \|f\|_{\mathcal H_x^{-2\alpha\eta-3\alpha+\frac{11}{2}}}.
\end{split}
\end{equation}
Since $0<\eta\leq \frac12$, we have $-2\alpha\eta-3\alpha+\frac{11}{2}\geq\frac{11}{2}-4\alpha$. It completes the proof of the lemma.

\cbdu

\medskip

\begin{Remark}
Recall $\mathcal H^{\frac72-2\alpha}$ is a critical space for (\ref{emhd}). We observe 
\[\frac{11}{2}-4\alpha<\frac72-2\alpha, \ \ \mbox{for}\ \ \alpha>1\]
 and $\frac{5}{2}-2\alpha<\frac{7}{2}-2\alpha$. Therefore, for $\max\{ \frac52-2\alpha, \frac{11}2-4\alpha\}\leq s<\frac72-2\alpha$, Lemmas \ref{le-heat3}, \ref{le-heat5} and \ref{le-heat6} hold for initial data $f$ in the supercritical regime. 
\end{Remark}

\medskip
 
Denote 
\begin{equation}\notag
\begin{split}
E_1(f, \beta, \alpha, \lambda)=&\ \left\{\omega\in\Omega: \|B_{f^\omega}\|_{ L^{\frac{1}{2\beta}}_t L^p_x}\geq \lambda \right\},\\
E_2(f, \beta, \alpha, \lambda)=&\ \left\{\omega\in\Omega: \|B_{f^\omega}\|_{ L^{(\beta, \frac{1}{\beta})}_t L^p_x}\geq \lambda \right\},\\
E_3(f, \beta, \alpha, \lambda)=&\ \left\{\omega\in\Omega: \|B_{f^\omega}\|_{L^{(\eta, 2)}_t \mathcal H^{\frac{11}{2}-2\alpha}_x}\geq \lambda \right\}.
\end{split}
\end{equation}
%We have the following probabilistic estimates.

\medskip

\begin{Lemma}\label{le-heat7}
Let $\alpha$, $\beta$, $\eta$ and $p$ satisfy the parameter conditions in Lemmas \ref{le-heat3}, \ref{le-heat5} and \ref{le-heat6}. Assume $f\in \mathcal H_x^{s}$ for $s\geq \max\{ \frac52-2\alpha, \frac{11}2-4\alpha\}$. There exist constants $c_1>0$ and $c_2>0$ such that   
\[P(E_i(f, \beta, \alpha, \lambda))\leq c_1e^{-\frac{c_2\lambda^2}{\|f\|^2_{\mathcal H_x^{s}}}},  \ \ \mbox{for}\ \ i=1, 2, 3.\]
\end{Lemma}
\pf
We only show the estimate for $E_1(f, \beta, \alpha, \lambda)$, since the other two can be handled analogously. In view of Bienaym\'e-Tchebishev's inequality, we deduce 
\begin{equation}\label{est-pp1}
P(E_1(f, \beta, \alpha, \lambda))= P\left(\left\{\omega\in\Omega: \|B_{f^\omega}\|_{ L^{\frac{1}{2\beta}}_t L^p_x}\geq \lambda \right\}\right)\leq \left(c_0\sqrt r\lambda^{-1}\|f\|_{\mathcal H_x^{s}}\right)^r
\end{equation}
for some constant $c_0>0$. If 
\[\left(\frac{\lambda}{c_0 e \|f\|_{\mathcal H_x^{s}}}\right)^2\geq \frac{1}{2\beta},\]
we take $r=\left(\frac{\lambda}{c_0 e \|f\|_{\mathcal H_x^{s}}}\right)^2$. It then follows from (\ref{est-pp1}) directly
\[P(E_i(f, \beta, \alpha, \lambda))\leq e^{-\frac{\lambda^2}{c_0^2 \|f\|^2_{\mathcal H_x^{s}}}}.\]
If otherwise 
\[\left(\frac{\lambda}{c_0 e \|f\|_{\mathcal H_x^{s}}}\right)^2< \frac{1}{2\beta}, \]
%\[\ \ i.e. \ \ c_0\lambda^{-1} \|f\|_{\mathcal H_x^{s}} \geq \sqrt{2\beta}e^{-1}\]
then there exists a constant $c_1$ such that we have from (\ref{est-pp1})
\[P(E_1(f, \beta, \alpha, \lambda))\leq c_1e^{-\frac{\lambda^2}{c_0^2 \|f\|^2_{\mathcal H_x^{s}}}},\]
by using the fact that $y^{-r}\leq c_1 e^{-y^2}$ for small $y>0$. 
\cbdu

\medskip

\medskip

%{\color{red} Lemma \ref{le-heat8} is not necessary.}
%\begin{Lemma}\label{le-heat8}
%Let $r\geq q\geq p\geq 2$ and $s\geq -\frac{2\alpha}{q}+\frac{1}{2}-\frac{1}{p}$. The free evolution $u_{g^\omega}$ satisfies
%\begin{equation}\notag
%\|u_{g^\omega}\|_{L^r_\omega L^q_t L^p_x}\lesssim \sqrt r \|g\|_{\mathcal H^s_x}.
%\end{equation}
%\end{Lemma}

%As a special case of Lemma \ref{le-heat2} and Lemma \ref{le-heat4} respectively, we also have
Note that the estimates in the lemmas above hold for $u_{g^\omega}$ with slight modifications. 
\begin{Lemma}\label{le-heat9}
Let $\frac{1}{2\alpha}+\frac{3}{2q\alpha}+2\gamma=1$. Assume $s> \frac32-2\alpha$. 
We have
\begin{equation}\notag
\begin{split}
\|u_{g^\omega}\|_{L^r_\omega L^{\frac{1}{2\gamma}}_t L^q_x}\lesssim&\ \sqrt r \|g\|_{\mathcal H^s_x}, \ \ \mbox{for} \ \ 
r\geq \frac{1}{2\gamma}\geq q\geq 2,\\
\|u_{g^\omega}\|_{L^r_\omega L^{(\gamma, \frac{1}{\gamma})}_t L^q_x}\lesssim&\ \sqrt r \|g\|_{\mathcal H^s_x}, \ \ \mbox{for} \ \ r\geq \frac{1}{\gamma}\geq q\geq 2.
\end{split}
\end{equation}
\end{Lemma}
%\begin{Lemma}\label{le-heat10}
%Let $r\geq \frac{1}{\gamma}\geq q\geq 2$ and $\frac{1}{2\alpha}+\frac{3}{2q\alpha}+2\gamma=1$. Let $s> \frac32-2\alpha$. We have the estimate
%\begin{equation}\notag
%\|u_{g^\omega}\|_{L^r_\omega L^{(\gamma, \frac{1}{\gamma})}_t L^q_x}\lesssim \sqrt r \|g\|_{\mathcal H^s_x}.
%\end{equation}
%\end{Lemma}

\medskip

\begin{Lemma}\label{le-heat11}
Let $r\geq4$ and $0<\zeta\leq\frac12$. Then
\begin{equation}\notag
\begin{split}
%\|u_{g^\omega}\|_{L^r_\omega L^{(\zeta, 2)}_t \mathcal H^{\frac{7}{2}-2\alpha}_x}\lesssim&\ \sqrt r \|g\|_{\mathcal H^s_x}, \ \ s\geq \frac{7}{2}-4\alpha,\\
\|u_{g^\omega}\|_{L^r_\omega L^{(\zeta, 4)}_t \mathcal H^{\frac{7}{2}-2\alpha}_x}\lesssim&\ \sqrt r \|g\|_{\mathcal H^s_x}, \ \ s\geq \frac{7}{2}-\frac{7\alpha}2.
\end{split}
\end{equation}
\end{Lemma}

Denote 
\begin{equation}\notag
\begin{split}
E_4(f, \beta, \alpha, \lambda)=&\ \left\{\omega\in\Omega: \|u_{g^\omega}\|_{ L^{\frac{1}{2\gamma}}_t L^q_x}\geq \lambda \right\},\\
E_5(f, \beta, \alpha, \lambda)=&\ \left\{\omega\in\Omega: \|u_{g^\omega}\|_{ L^{(\gamma, \frac{1}{\gamma})}_t L^q_x}\geq \lambda \right\},\\
E_6(f, \beta, \alpha, \lambda)=&\ \left\{\omega\in\Omega: \|u_{g^\omega}\|_{L^{(\zeta, 4)}_t \mathcal H^{\frac{7}{2}-2\alpha}_x}\geq \lambda \right\}.
\end{split}
\end{equation}

\begin{Lemma}\label{le-heat12}
Let $\alpha$, $\gamma$, $\zeta$ and $q$ satisfy the parameter conditions in Lemmas \ref{le-heat9} and \ref{le-heat11}. Assume $f\in \mathcal H_x^{s}$ for $s\geq\max\{\frac{7}{2}-\frac72\alpha, \frac{3}{2}-2\alpha\}$. There exist constants $c_3>0$ and $c_4>0$ such that   
\[P(E_i(f, \beta, \alpha, \lambda))\leq c_3e^{-\frac{c_4\lambda^2}{\|g\|^2_{\mathcal H_x^{s}}}},  \ \ \mbox{for}\ \ i=4, 5, 6.\]
\end{Lemma}

\bigskip

\section{Well-posedness of the electron MHD}\label{sec-emhd}

We prove Theorem \ref{thm1} for the electron MHD (\ref{emhd}) in this section. As discussed earlier, the electron MHD contains the nonlinear term with the highest derivative from the Hall MHD system. Hence we encounter the most challenging estimates in this part. 

To take the advantage of the improved $L^p$ estimates for the free evolution $B_{f^\omega}=e^{-t(-\Delta)^\alpha} f^\omega$, we look for a solution of (\ref{emhd}) with initial data $f^\omega$ in the form 
\[B=B_{f^\omega}+H\] with the nonlinear part $H$ solving the Cauchy problem
%Denote the bilinear operator
%\begin{equation}\notag
%\mathcal B(u,v)=(\nabla\times u)\times v.
%\end{equation}
%Note that  \[\mathcal B(u,u)= (u\cdot\nabla )u-\nabla \frac{|u|^2}{2}. \]
%If $\nabla\cdot u=0$, we can further write \[\mathcal B(u,u)= \nabla\cdot(u\otimes u)-\nabla \frac{|u|^2}{2}. \]
%One can check that if $B$ satisfies (\ref{emhd}) with initial data $f^\omega$, the nonlinear part $H$ solves the Cauchy problem
\begin{equation}\label{H}
\begin{split}
H_t+\nabla\times\nabla\cdot \left((B_{f^\omega}+H)\otimes (B_{f^\omega}+H)\right)
=&-(-\Delta)^\alpha H,\\
\nabla\cdot H=&\ 0,\\
H(x, 0)=&\ 0.
\end{split}
\end{equation}
Here we used the rewriting 
\[\nabla\times[(\nabla\times (B_{f^\omega}+H))\times (B_{f^\omega}+H)]= \nabla\times\nabla\cdot \left((B_{f^\omega}+H)\otimes (B_{f^\omega}+H)\right)\]
since $\nabla \cdot H=0$ and $\nabla\cdot B_{f^\omega}=0$.

Obviously Theorem \ref{thm1} follows from the well-posedness of (\ref{H}) in a suitable subspace of $C([0,T]; \mathcal H^{\frac{7}{2}-2\alpha}(\mathbb T^3))$. Therefore we only need to show:

%\medskip

\begin{Theorem}\label{thm-H}
Let $f\in \mathcal H^{s}$ with $s\geq\max\{\frac{11}{2}-4\alpha, \frac{5}{2}-2\alpha\}$ and $f$ be zero-mean. Let $\beta>0$ and $p\geq 2$ satisfy $\frac{1}{\alpha}+\frac{3}{2p\alpha}+2\beta=1$. There exists a set $\Sigma\subset\Omega$ with $P(\Sigma)=1$ such that for any $\omega\in\Sigma$ system (\ref{H}) has a unique solution $H$ satisfying 
\[H\in C([0,T]; \mathcal H^{\frac{7}{2}-2\alpha}(\mathbb T^3))\cap L^{(\beta,\frac{1}{\beta})}(0, T; L^p(\mathbb T^3))\cap L^{\frac{1}{2\beta}}(0, T; L^p(\mathbb T^3)) \] 
for some $T>0$. If in addition, $\|f\|_{\mathcal H^{s}}\leq c$ for some constant $c>0$, then $T=\infty$. 
\end{Theorem}
\pf
We proceed by employing a fixed point argument. 
Denote 
\begin{equation}\notag
Q(x,t)= (B_{f^\omega}+H)\otimes (B_{f^\omega}+H) (x, t).
\end{equation}
The integral form of (\ref{H}) is given by
\begin{equation}\label{H2}
\begin{split}
H(x,t)
%=&-\int_0^t e^{-(t-s)(-\Delta)^\alpha} \tilde Q(x, s)\, ds\\
=&-\int_0^t e^{-(t-\tau)(-\Delta)^\alpha} \nabla\times\nabla\cdot Q(x, \tau)\, d\tau.
\end{split}
\end{equation}
Denote the map 
\begin{equation}\notag
\Phi(H)(t)=-\int_0^t e^{-(t-\tau)(-\Delta)^\alpha} \nabla\times\nabla\cdot Q(x, \tau)\, d\tau.
\end{equation}
Define the subspace $\mathcal Y\subset C([0,T]; \mathcal H^{\frac{7}{2}-2\alpha}(\mathbb T^3))$ as
\begin{equation}\notag
\mathcal Y=C([0,T]; \mathcal H^{\frac{7}{2}-2\alpha}(\mathbb T^3))\cap L^{(\beta,\frac{1}{\beta})}(0, T; L^p(\mathbb T^3))\cap L^{\frac{1}{2\beta}}(0, T; L^p(\mathbb T^3)).
\end{equation}
We claim that the map $\Phi$ is a contraction on $\mathcal Y$ by showing that:\\
(i) $\Phi$ maps $\mathcal Y$ onto itself;\\
(ii) For any $H_1\in\mathcal Y$ and $H_2\in \mathcal Y$, we have
\[\|\Phi(H_1)-\Phi(H_2)\|_{\mathcal Y}\leq c \|H_1-H_2\|_{\mathcal Y}.\]
In order to show (i), we estimate $\|\Phi(H)\|_{\mathcal H^{\frac{7}{2}-2\alpha}(\mathbb T^3)}$, $\|\Phi(H)\|_{L^{(\beta,\frac{1}{\beta})}(0, T; L^p(\mathbb T^3))}$ and $\|\Phi(H)\|_{L^{\frac{1}{2\beta}}(0, T; L^p(\mathbb T^3))}$ respectively in the following. We first expand $\Phi(H)$ as 
\begin{equation}\notag
\begin{split}
\Phi(H)(t)
=&-\int_0^t e^{-(t-\tau)(-\Delta)^\alpha} \nabla\times\nabla\cdot (H(x, \tau)\otimes H(x,\tau))\, d\tau\\
&-\int_0^t e^{-(t-\tau)(-\Delta)^\alpha} \nabla\times\nabla\cdot (H(x, \tau)\otimes B_{f^\omega}(x,\tau))\, d\tau\\
&-\int_0^t e^{-(t-\tau)(-\Delta)^\alpha} \nabla\times\nabla\cdot (B_{f^\omega}(x, \tau)\otimes H(x,\tau))\, d\tau\\
&-\int_0^t e^{-(t-\tau)(-\Delta)^\alpha} \nabla\times\nabla\cdot (B_{f^\omega}(x, \tau)\otimes B_{f^\omega}(x,\tau))\, d\tau\\
=:& -\Phi_1-\Phi_2-\Phi_3-\Phi_4.
\end{split}
\end{equation}

We estimate  $\Phi_1$ in $\mathcal H^{\frac{7}{2}-2\alpha}(\mathbb T^3)$ as
\begin{equation}\notag
\begin{split}
&\|\Phi_1\|_{\mathcal H^{\frac{7}{2}-2\alpha}(\mathbb T^3)}\\
=&\ \left\| (1-\Delta)^{\frac74-\alpha}\int_0^t e^{-(t-\tau)(-\Delta)^\alpha} \nabla\times\nabla\cdot (H(x, \tau)\otimes H(x,\tau))\, d\tau \right\|_{L^2_x}\\
=&\ \left\| \int_0^t e^{-(t-\tau)(-\Delta)^\alpha}(-\Delta) (1-\Delta)^{\frac74-\alpha}\frac{\nabla\times\nabla\cdot}{(-\Delta)} (H(x, \tau)\otimes H(x,\tau))\, d\tau \right\|_{L^2_x}\\
\lesssim&\  \int_0^t \left\|e^{-(t-\tau)(-\Delta)^\alpha}(-\Delta) (1-\Delta)^{\frac74-\alpha} (H(x, \tau)\otimes H(x,\tau)) \right\|_{L^2_x}\, d\tau \\
\lesssim&\  \int_0^t \left\|e^{-(t-\tau)|\xi|^{2\alpha}}|\xi|^{2} \mathcal F\left((1-\Delta)^{\frac74-\alpha} (H(x, \tau)\otimes H(x,\tau))\right) \right\|_{L^2_{\xi}}\, d\tau \\
\end{split}
\end{equation}
where Plancherel's theorem was applied in the last step. Using H\"older's inequality and Lemma \ref{le-heat1} we obtain
\begin{equation}\notag
\begin{split}
&\left\|e^{-(t-\tau)|\xi|^{2\alpha}}|\xi|^{2} \mathcal F\left((1-\Delta)^{\frac74-\alpha} (H(x, \tau)\otimes H(x,\tau))\right) \right\|_{L^2_{\xi}}\\
\lesssim&\ \left\|e^{-(t-\tau)|\xi|^{2\alpha}}|\xi|^{2}\right\|_{L^p_{\xi}} \left\|\mathcal F\left((1-\Delta)^{\frac74-\alpha} (H(x, \tau)\otimes H(x,\tau))\right) \right\|_{L^{\frac{2p}{p-2}}_{\xi}}\\
\lesssim&\ (t-\tau)^{-\frac{1}{\alpha}-\frac{3}{2p\alpha}} \left\|(1-\Delta)^{\frac74-\alpha} (H(x, \tau)\otimes H(x,\tau)) \right\|_{L^{\frac{2p}{p+2}}_{x}}\\
\lesssim&\ (t-\tau)^{-\frac{1}{\alpha}-\frac{3}{2p\alpha}} \left\| H\nabla^{\frac72-2\alpha}H \right\|_{L^{\frac{2p}{p+2}}_{x}}\\
\lesssim&\ (t-\tau)^{-\frac{1}{\alpha}-\frac{3}{2p\alpha}} \left\| H \right\|_{L^{p}_{x}} \left\| \nabla^{\frac72-2\alpha}H \right\|_{L^{2}_{x}}.
\end{split}
\end{equation}
It follows from the last two inequalities that
\begin{equation}\notag
\begin{split}
&\|\Phi_1\|_{\mathcal H^{\frac{7}{2}-2\alpha}(\mathbb T^3)}\\
\lesssim&\ \int_0^t (t-\tau)^{-\frac{1}{\alpha}-\frac{3}{2p\alpha}} \left\| H \right\|_{L^{p}_{x}} \left\| \nabla^{\frac72-2\alpha}H \right\|_{L^{2}_{x}}\, d\tau\\
\lesssim&\ \left\|H \right\|_{L^\infty_t\mathcal H^{\frac72-2\alpha}_{x}}\int_0^t (t-\tau)^{-\frac{1}{\alpha}-\frac{3}{2p\alpha}}\tau^{-\beta} \left(\tau^\beta\left\| H \right\|_{L^{p}_{x}}\right) \, d\tau\\
\lesssim&\ \left\|H \right\|_{L^\infty_t\mathcal H^{\frac72-2\alpha}_{x}}\left(\int_0^t (t-\tau)^{\left(-\frac{1}{\alpha}-\frac{3}{2p\alpha}\right)\frac{1}{1-\beta}}\tau^{-\frac{\beta}{1-\beta}} \, d\tau\right)^{1-\beta} \left(\int_0^t \tau \| H \|_{L^{p}_{x}}^{\frac1\beta}\, d\tau\right)^{\beta}.
\end{split}
\end{equation}
Based on the assumptions on the parameters, we observe that 
\[\left(\frac{1}{\alpha}+\frac{3}{2p\alpha}\right)\frac{1}{1-\beta}+\frac{\beta}{1-\beta}=1,\]
and 
\[0<\left(\frac{1}{\alpha}+\frac{3}{2p\alpha}\right)\frac{1}{1-\beta}<1, \ \ 0<\frac{\beta}{1-\beta}<1.\]
Hence Lemma \ref{le-beta} implies the time integral is bounded.
Therefore we conclude
\begin{equation}\label{est-phi1}
\|\Phi_1\|_{\mathcal H^{\frac{7}{2}-2\alpha}(\mathbb T^3)}
\lesssim \left\|H \right\|_{L^\infty_t\mathcal H^{\frac72-2\alpha}_{x}}\left\|H \right\|_{L^{(\beta, \frac{1}{\beta})}_t L^p_{x}}.
\end{equation}

We continue to estimate $\Phi_1$ in $L^{(\beta, \frac{1}{\beta})}_t L^p_{x}$,
\begin{equation}\notag
\begin{split}
&\|\Phi_1\|_{L^{(\beta, \frac{1}{\beta})}_t L^p_{x}}\\
=&\ \left\|t^\beta \int_0^t e^{-(t-\tau)(-\Delta)^\alpha} \nabla\times\nabla\cdot (H(x, \tau)\otimes H(x,\tau))\, d\tau \right\|_{L^{\frac{1}{\beta}}_t L^p_{x}}\\
\lesssim &\ \left\|t^\beta \int_0^t \left\|e^{-(t-\tau)(-\Delta)^\alpha} (-\Delta) (H(x, \tau)\otimes H(x,\tau))\right\|_{L^p_{x}}\, d\tau \right\|_{L^{\frac{1}{\beta}}_t}.
\end{split}
\end{equation}
Applying H\"older's inequality and Lemma \ref{le-heat1} yields
\begin{equation}\notag
\begin{split}
&\left\|e^{-(t-\tau)(-\Delta)^\alpha} (-\Delta) (H(x, \tau)\otimes H(x,\tau))\right\|_{L^p_{x}}\\
\lesssim&\ \left\|e^{-(t-\tau)|\xi|^{2\alpha}} |\xi|^2 \mathcal F\left((H(x, \tau)\otimes H(x,\tau))\right)\right\|_{L^{\frac{p}{p-1}}_{\xi}}\\
\lesssim&\ \left\|e^{-(t-\tau)|\xi|^{2\alpha}} |\xi|^2\right\|_{L^{p}_{\xi}}\left\| \mathcal F\left((H(x, \tau)\otimes H(x,\tau))\right)\right\|_{L^{\frac{p}{p-2}}_{\xi}}\\
\lesssim&\ (t-\tau)^{-\frac1{\alpha}-\frac{3}{2p\alpha}}\left\|H(x, \tau)\otimes H(x,\tau)\right\|_{L^{\frac{p}{2}}_{x}}\\
\lesssim&\ (t-\tau)^{-\frac1{\alpha}-\frac{3}{2p\alpha}}\left\|H\right\|^2_{L^{p}_{x}}.
\end{split}
\end{equation}
Combining the last two inequalities we infer
\begin{equation}\notag
\begin{split}
&\|\Phi_1\|_{L^{(\beta, \frac{1}{\beta})}_t L^p_{x}}\\
\lesssim &\ \left\|\int_0^t t^\beta(t-\tau)^{-\frac1{\alpha}-\frac{3}{2p\alpha}}\left\|H\right\|^2_{L^{p}_{x}} \, d\tau \right\|_{L^{\frac{1}{\beta}}_t}\\
\lesssim &\ \left\| \int_0^{\frac{t}{2}} t^\beta (t-\tau)^{-\frac1{\alpha}-\frac{3}{2p\alpha}}\left\|H\right\|^2_{L^{p}_{x}} \, d\tau \right\|_{L^{\frac{1}{\beta}}_t}\\
&+ \left\|\int_{\frac{t}{2}}^t t^\beta(t-\tau)^{-\frac1{\alpha}-\frac{3}{2p\alpha}}\left\|H\right\|^2_{L^{p}_{x}} \, d\tau \right\|_{L^{\frac{1}{\beta}}_t}\\
=:&\ I_1+I_2.
\end{split}
\end{equation}
The term $I_1$ can be estimated by using Young's inequality
\begin{equation}\notag
\begin{split}
& \left\| \int_0^{\frac{t}{2}} t^\beta (t-\tau)^{-\frac1{\alpha}-\frac{3}{2p\alpha}}\left\|H\right\|^2_{L^{p}_{x}} \, d\tau \right\|_{L^{\frac{1}{\beta}}_t}\\
\lesssim&\ \left\| \int_0^{\frac{t}{2}}  (t-\tau)^{\beta-\frac1{\alpha}-\frac{3}{2p\alpha}}\left\|H\right\|^2_{L^{p}_{x}} \, d\tau \right\|_{L^{\frac{1}{\beta}}_t}\\
\lesssim&\ \left(\int_0^{\frac{t}{2}} (t-\tau)^{(\beta-\frac{1}{\alpha}-\frac{3}{2p\alpha})(\frac{1}{1-3\beta})}\, d\tau\right)^{1-3\beta}
\left(\int_0^{\frac{t}{2}}\left\|H\right\|^{2\cdot \frac{1}{4\beta}}_{L^{p}_{x}}\, d\tau \right)^{4\beta}\\
\lesssim&\ \left(\int_0^{\frac{t}{2}} (t-\tau)^{-1}\, d\tau\right)^{1-3\beta} \|H\|^2_{L^{\frac{1}{2\beta}}_t L^{p}_{x}}\\
\lesssim&\ \|H\|^2_{L^{\frac{1}{2\beta}}_t L^{p}_{x}}.
\end{split}
\end{equation}
While the term $I_2$ is estimated by employing H\"older's inequality and Hardy-Littlewood-Sobolev lemma
\begin{equation}\notag
\begin{split}
& \left\| \int_{\frac{t}{2}}^t t^\beta (t-\tau)^{-\frac1{\alpha}-\frac{3}{2p\alpha}}\left\|H\right\|^2_{L^{p}_{x}} \, d\tau \right\|_{L^{\frac{1}{\beta}}_t}\\
\lesssim&\ \left\| \int_0^{\frac{t}{2}}  (t-\tau)^{-\frac1{\alpha}-\frac{3}{2p\alpha}}\tau^\beta \|H\|_{L^{p}_{x}} \|H\|_{L^{p}_{x}}\, d\tau \right\|_{L^{\frac{1}{\beta}}_t}\\
\lesssim&\ \left\| \left(\int_{\frac{t}{2}}^t\tau \|H\|_{L^{p}_{x}}^{\frac{1}{\beta}} \, d\tau\right)^{\beta}\left(\int_{\frac{t}{2}}^t  (t-\tau)^{(-\frac1{\alpha}-\frac{3}{2p\alpha})\frac{1}{1-\beta}} \|H\|_{L^{p}_{x}}^{\frac{1}{1-\beta}}\, d\tau\right)^{1-\beta} \right\|_{L^{\frac{1}{\beta}}_t}\\
\lesssim&\ \|H\|_{L^{(\beta, \frac{1}{\beta})}_t L^{p}_{x}}\left\| \int_{\frac{t}{2}}^t  (t-\tau)^{(-\frac1{\alpha}-\frac{3}{2p\alpha})\frac{1}{1-\beta}} \|H\|_{L^{p}_{x}}^{\frac{1}{1-\beta}}\, d\tau \right\|^{1-\beta}_{L^{\frac{1-\beta}{\beta}}_t}\\
\lesssim&\  \|H\|_{L^{(\beta, \frac{1}{\beta})}_t L^{p}_{x}} 
\left\| \|H\|_{L^{p}_{x}}^{\frac{1}{1-\beta}} \right\|_{L^{\frac{1-\beta}{2\beta}}_t}^{1-\beta}\\
%\lesssim&\  \|H\|_{L^{(\beta, \frac{1}{\beta})}_t L^{p}_{x}} \left(\int_{\frac{t}{2}}^t\left\|H\right\|^{(\frac{1}{1-\beta})(\frac{1-\beta}{2\beta})}_{L^{p}_{x}}\, d\tau \right)^{2\beta}\\
\lesssim&\ \|H\|_{L^{(\beta, \frac{1}{\beta})}_t L^{p}_{x}}\|H\|_{L^{\frac{1}{2\beta}}_t L^{p}_{x}}.
\end{split}
\end{equation}
Putting together the estimates above we get
\begin{equation}\label{est-phi2}
\|\Phi_1\|_{L^{(\beta, \frac{1}{\beta})}_t L^p_{x}}
\lesssim \left\|H \right\|_{L^{\frac{1}{2\beta}}_tL^p_x}^2+\left\|H \right\|_{L^{(\beta, \frac{1}{\beta})}_t L^p_{x}}\left\|H \right\|_{L^{\frac{1}{2\beta}}_tL^p_x}.
\end{equation}

\medskip

While to estimate $\Phi_1$ in $L^{\frac{1}{2\beta}}_t L^p_{x}$, we start with
\begin{equation}\notag
\begin{split}
&\|\Phi_1\|_{L^{\frac{1}{2\beta}}_t L^p_{x}}\\
=&\ \left\| \int_0^t e^{-(t-\tau)(-\Delta)^\alpha} \nabla\times\nabla\cdot (H(x, \tau)\otimes H(x,\tau))\, d\tau \right\|_{L^{\frac{1}{2\beta}}_t L^p_{x}}\\
\lesssim &\ \left\| \int_0^t \left\|e^{-(t-\tau)(-\Delta)^\alpha} (-\Delta) (H(x, \tau)\otimes H(x,\tau))\right\|_{L^p_{x}}\, d\tau \right\|_{L^{\frac{1}{2\beta}}_t}.
\end{split}
\end{equation}
Similarly as before, we deduce by applying H\"older's inequality and Lemma \ref{le-heat1} 
\begin{equation}\notag
\begin{split}
& \left\|e^{-(t-\tau)(-\Delta)^\alpha} (-\Delta) (H(x, \tau)\otimes H(x,\tau))\right\|_{L^p_{x}}\\
\lesssim &\ \left\|e^{-(t-\tau)|\xi|^{2\alpha}} |\xi|^2 \mathcal F\left(H(x, \tau)\otimes H(x,\tau)\right)\right\|_{L^{\frac{p}{p-1}}_{\xi}}\\
\lesssim &\ \left\|e^{-(t-\tau)|\xi|^{2\alpha}} |\xi|^2\right\|_{L^{p}_{\xi}} \left\|\mathcal F\left(H(x, \tau)\otimes H(x,\tau)\right)\right\|_{L^{\frac{p}{p-2}}_{\xi}}\\
\lesssim &\ (t-\tau)^{-\frac{1}{\alpha}-\frac{3}{2p\alpha}} \left\|H(x, \tau)\otimes H(x,\tau)\right\|_{L^{\frac{p}{2}}_{x}}\\
\lesssim &\ (t-\tau)^{-\frac{1}{\alpha}-\frac{3}{2p\alpha}} \|H\|^2_{L^{p}_{x}}.
\end{split}
\end{equation}
Invoking Hardy-Littlewood-Sobolev lemma again, the last two inequalities together imply
\begin{equation}\label{est-phi3}
\|\Phi_1\|_{L^{\frac{1}{2\beta}}_t L^p_{x}}
\lesssim \left\| \int_0^t (t-\tau)^{-\frac{1}{\alpha}-\frac{3}{2p\alpha}} \|H\|^2_{L^{p}_{x}} \, d\tau \right\|_{L^{\frac{1}{2\beta}}_t}
\lesssim  \left\|H \right\|_{L^{\frac{1}{2\beta}}_tL^p_x}^2.
\end{equation}

The estimate of $\Phi_4$ is analogous to that of $\Phi_1$. We sketch some details below. 
Regarding the estimate in $\mathcal H^{\frac{7}{2}-2\alpha}(\mathbb T^3)$, we split $\Phi_4$ as
\begin{equation}\notag
\begin{split}
&\|\Phi_4\|_{\mathcal H^{\frac{7}{2}-2\alpha}(\mathbb T^3)}\\
=&\ \left\| (1-\Delta)^{\frac74-\alpha}\int_0^t e^{-(t-\tau)(-\Delta)^\alpha} \nabla\times\nabla\cdot (B_{f^\omega}(x, \tau)\otimes B_{f^\omega}(x,\tau))\, d\tau \right\|_{L^2_x}\\
%=&\ \left\| \int_0^t e^{-(t-\tau)(-\Delta)^\alpha}(-\Delta) (1-\Delta)^{\frac74-\alpha}\frac{\nabla\times\nabla\cdot}{(-\Delta)} (B_{f^\omega}(x, \tau)\otimes B_{f^\omega}(x,\tau))\, d\tau \right\|_{L^2_x}\\
%\lesssim&\  \int_0^t \left\|e^{-(t-\tau)(-\Delta)^\alpha}(-\Delta) (1-\Delta)^{\frac74-\alpha} (B_{f^\omega}(x, \tau)\otimes B_{f^\omega}(x,\tau)) \right\|_{L^2_x}\, d\tau \\
\lesssim&\  \int_0^\delta \left\|e^{-(t-\tau)|\xi|^{2\alpha}}|\xi|^{2+\frac72-2\alpha} \mathcal F\left(B_{f^\omega}(x, \tau)\otimes B_{f^\omega}(x,\tau)\right) \right\|_{L^2_{\xi}}\, d\tau \\
&+ \int_\delta^t \left\|e^{-(t-\tau)|\xi|^{2\alpha}} \mathcal F\left((-\Delta)(1-\Delta)^{\frac74-\alpha} (B_{f^\omega}(x, \tau)\otimes B_{f^\omega}(x,\tau))\right) \right\|_{L^2_{\xi}}\, d\tau \\
=:&\ I_3+I_4
\end{split}
\end{equation}
where $0<\delta<t$ is a small constant. Applying Lemma \ref{le-heat1} we infer for $2\leq p\leq 4$
\begin{equation}\notag
\begin{split}
&\left\|e^{-(t-\tau)|\xi|^{2\alpha}}|\xi|^{2+\frac72-2\alpha} \mathcal F\left(B_{f^\omega}(x, \tau)\otimes B_{f^\omega}(x,\tau)\right) \right\|_{L^2_{\xi}}\\
\lesssim &\ \left\|e^{-(t-\tau)|\xi|^{2\alpha}}|\xi|^{2+\frac72-2\alpha}\right\|_{L^{\frac{2p}{4-p}}_{\xi}}\left\|\mathcal F\left(B_{f^\omega}(x, \tau)\otimes B_{f^\omega}(x,\tau)\right) \right\|_{L^{\frac{p}{p-2}}_{\xi}}\\
\lesssim &\ (t-\tau)^{-\frac{11}{4\alpha}+1-\frac{3}{2\alpha}\cdot\frac{4-p}{2p}}\left\|B_{f^\omega}(x, \tau)\otimes B_{f^\omega}(x,\tau) \right\|_{L^{\frac{p}{2}}_{x}}\\
\lesssim &\ (t-\tau)^{-\frac{11}{4\alpha}+1-\frac{3(4-p)}{4p\alpha}}\|B_{f^\omega} \|^2_{L^{p}_{x}}.
\end{split}
\end{equation}
Thus, applying H\"older's inequality and Lemma \ref{le-beta} leads to
\begin{equation}\notag
\begin{split}
I_3\lesssim&\ \int_0^\delta (t-\tau)^{-\frac{11}{4\alpha}+1-\frac{3(4-p)}{4p\alpha}}\|B_{f^\omega} \|^2_{L^{p}_{x}}\, d\tau\\
=&\ \int_0^\delta \left((t-\tau)^{-\frac{11}{4\alpha}+1-\frac{3(4-p)}{4p\alpha}}\tau^{-2\beta}\right) \left(\tau^{\beta}\|B_{f^\omega} \|_{L^{p}_{x}}\right)^2\, d\tau\\
\lesssim&\ \left(\int_0^\delta (t-\tau)^{(-\frac{11}{4\alpha}+1-\frac{3(4-p)}{4p\alpha})\frac{1}{1-2\beta}}\tau^{-2\beta\cdot\frac{1}{1-2\beta}} \, d\tau\right)^{1-2\beta} \|B_{f^\omega}\|^2_{L^{(\beta,\frac{1}{\beta})}_tL^p_x}\\
=&\ \left(\int_0^\delta (t-\tau)^{-1+\frac{2\beta}{1-2\beta}}\tau^{-\frac{2\beta}{1-2\beta}} \, d\tau\right)^{1-2\beta} \|B_{f^\omega}\|^2_{L^{(\beta,\frac{1}{\beta})}_tL^p_x}\\
\lesssim&\ \|B_{f^\omega}\|^2_{L^{(\beta,\frac{1}{\beta})}_tL^p_x}.
\end{split}
\end{equation}

On the other hand, for $p>4$ we have
\begin{equation}\notag
\begin{split}
&\left\|e^{-(t-\tau)|\xi|^{2\alpha}}|\xi|^{2+\frac72-2\alpha} \mathcal F\left(B_{f^\omega}(x, \tau)\otimes B_{f^\omega}(x,\tau)\right) \right\|_{L^2_{\xi}}\\
\lesssim &\ \left\|e^{-(t-\tau)|\xi|^{2\alpha}}|\xi|^{2+\frac72-2\alpha}\right\|_{L^{\infty}_{\xi}}\left\|\mathcal F\left(B_{f^\omega}(x, \tau)\otimes B_{f^\omega}(x,\tau)\right) \right\|_{L^{2}_{\xi}}\\
\lesssim &\ (t-\tau)^{-\frac{11}{4\alpha}+1}\left\|B_{f^\omega}(x, \tau)\otimes B_{f^\omega}(x,\tau) \right\|_{L^{2}_{x}}\\
\lesssim &\ (t-\tau)^{-\frac{11}{4\alpha}+1}\|B_{f^\omega} \|^2_{L^{4}_{x}}
\end{split}
\end{equation}
and hence
\begin{equation}\notag
\begin{split}
I_3\lesssim&\ \int_0^\delta (t-\tau)^{-\frac{11}{4\alpha}+1}\|B_{f^\omega} \|^2_{L^{4}_{x}}\, d\tau\\
=&\ \int_0^\delta \left((t-\tau)^{-\frac{11}{4\alpha}+1}\tau^{-2\beta}\right) \left(\tau^{\beta}\|B_{f^\omega} \|_{L^{4}_{x}}\right)^2\, d\tau\\
\lesssim&\ \left(\int_0^\delta (t-\tau)^{(-\frac{11}{4\alpha}+1)\frac{1}{1-2\beta}}\tau^{-\frac{2\beta}{1-2\beta}} \, d\tau\right)^{1-2\beta} \|B_{f^\omega}\|^2_{L^{(\beta,\frac{1}{\beta})}_tL^4_x}\\
\lesssim&\ \|B_{f^\omega}\|^2_{L^{(\beta,\frac{1}{\beta})}_tL^4_x}.
\end{split}
\end{equation}
In the estimate above, the time integral is handled as
\begin{equation}\notag
\begin{split}
&\int_0^\delta (t-\tau)^{(-\frac{11}{4\alpha}+1)\frac{1}{1-2\beta}}\tau^{-\frac{2\beta}{1-2\beta}} \, d\tau\\
< &\ t^{(-\frac{11}{4\alpha}+1)\frac{1}{1-2\beta}} \int_0^\delta \tau^{-\frac{2\beta}{1-2\beta}} \, d\tau\\
\lesssim&\ 1
\end{split}
\end{equation}
since $t>\delta>0$ and observing that for $1<\alpha<\frac74$
\[(-\frac{11}{4\alpha}+1)\frac{1}{1-2\beta}<0,\ \ \  0<\frac{2\beta}{1-2\beta}<1.\]

We continue to estimate $I_4$. Again invoking H\"older's inequality and Lemma \ref{le-heat1} we obtain
\begin{equation}\notag
\begin{split}
&\left\|e^{-(t-\tau)|\xi|^{2\alpha}} \mathcal F\left((-\Delta)(1-\Delta)^{\frac74-\alpha} (B_{f^\omega}(x, \tau)\otimes B_{f^\omega}(x,\tau))\right) \right\|_{L^2_{\xi}}\\
\lesssim&\ \left\|e^{-(t-\tau)|\xi|^{2\alpha}} \right\|_{L^{p}_{\xi}} \left\|\mathcal F\left((-\Delta)(1-\Delta)^{\frac74-\alpha} (B_{f^\omega}(x, \tau)\otimes B_{f^\omega}(x,\tau))\right) \right\|_{L^{\frac{2p}{p-2}}_{\xi}}\\
\lesssim&\ (t-\tau)^{-\frac{3}{2p\alpha}} \left\|(-\Delta)(1-\Delta)^{\frac74-\alpha} (B_{f^\omega}(x, \tau)\otimes B_{f^\omega}(x,\tau)) \right\|_{L^{\frac{2p}{p+2}}_{x}}\\
\lesssim&\ (t-\tau)^{-\frac{3}{2p\alpha}} \|\nabla^{\frac{11}{2}-2\alpha} B_{f^\omega}\|_{L^{2}_{x}}  \|B_{f^\omega}\|_{L^{p}_{x}}.
\end{split}
\end{equation}
Therefore we have for $\frac1{\alpha}-\frac12<\eta<\frac12-2\beta$ and $1<\alpha<\frac74$
\begin{equation}\notag
\begin{split}
I_4\lesssim&\ \int_\delta^t (t-\tau)^{-\frac{3}{2p\alpha}} \|\nabla^{\frac{11}{2}-2\alpha} B_{f^\omega}\|_{L^{2}_{x}}  \|B_{f^\omega}\|_{L^{p}_{x}} \, d\tau\\
=&\ \int_\delta^t (t-\tau)^{-\frac{3}{2p\alpha}}\tau^{-\eta-\beta} \left(\tau^\eta\|\nabla^{\frac{11}{2}-2\alpha} B_{f^\omega}\|_{L^{2}_{x}}\right)\left(\tau^{\beta}  \|B_{f^\omega}\|_{L^{p}_{x}}\right) \, d\tau\\
\lesssim &\ \left(\int_\delta^t (t-\tau)^{-\frac{3}{2p\alpha}\cdot \frac{2}{1-2\beta}}\tau^{(-\eta-\beta)\frac{2}{1-2\beta}} \, d\tau\right)^{\frac{1-2\beta}{2}}\|B_{f^\omega}\|_{L^{(\eta, 2)}_t \mathcal H^{\frac{11}{2}-2\alpha}_x} \|B_{f^\omega}\|_{L^{(\beta, \frac{1}{\beta})}_t L^{p}_x}\\
\lesssim &\ \|B_{f^\omega}\|_{L^{(\eta, 2)}_t \mathcal H^{\frac{11}{2}-2\alpha}_x} \|B_{f^\omega}\|_{L^{(\beta, \frac{1}{\beta})}_t L^{p}_x}\\
\end{split}
\end{equation}
where the time integral is estimated as
\begin{equation}\notag
\begin{split}
&\int_\delta^t (t-\tau)^{-\frac{3}{2p\alpha}\cdot \frac{2}{1-2\beta}}\tau^{(-\eta-\beta)\frac{2}{1-2\beta}} \, d\tau\\
=&\ t^{-\frac{\frac12+\eta-\frac1{\alpha}}{1-2\beta}}\int_{\delta/t}^{1} (1-\tau')^{- \frac{3}{p\alpha(1-2\beta)}} (\tau')^{-\frac{2(\eta+\beta)}{1-2\beta}} \, d\tau', \\
\lesssim &\ B\left(1-\frac{2(\eta+\beta)}{1-2\beta}, 1- \frac{3}{p\alpha(1-2\beta)}\right)\\
\lesssim &\ 1 
\end{split}
\end{equation}
for $t\geq \delta>0$, since
\[\frac{\frac12+\eta-\frac1{\alpha}}{1-2\beta}>0, \ \ 0<\frac{2(\eta+\beta)}{1-2\beta}<1, \ \ 0<\frac{3}{p\alpha(1-2\beta)}<1.\]
Summarizing the estimates above leads to
\begin{equation}\label{est-phi4}
\begin{split}
\|\Phi_4\|_{\mathcal H^{\frac{7}{2}-2\alpha}(\mathbb T^3)}\lesssim&\ \|B_{f^\omega}\|^2_{L^{(\beta,\frac{1}{\beta})}_tL^p_x}
+\|B_{f^\omega}\|^2_{L^{(\beta,\frac{1}{\beta})}_tL^4_x}\\
&+\|B_{f^\omega}\|_{L^{(\eta, 2)}_t \mathcal H^{\frac{11}{2}-2\alpha}_x} \|B_{f^\omega}\|_{L^{(\beta, \frac{1}{\beta})}_t L^{p}_x}.
\end{split}
\end{equation}
The estimates of $\Phi_4$ in $L^{(\beta, \frac{1}{\beta})}_t L^{p}_x$ and $L^{\frac{1}{2\beta}}_t L^{p}_x$ are identical to that of $\Phi_1$. Hence we collect the estimates
\begin{equation}\label{est-phi5}
\begin{split}
\|\Phi_4\|_{L^{(\beta, \frac{1}{\beta})}_t L^p_{x}}
\lesssim&\ \left\|B_{f^\omega} \right\|_{L^{\frac{1}{2\beta}}_tL^p_x}^2+\left\|B_{f^\omega} \right\|_{L^{(\beta, \frac{1}{\beta})}_t L^p_{x}}\left\|B_{f^\omega} \right\|_{L^{\frac{1}{2\beta}}_tL^p_x},\\
\|\Phi_4\|_{L^{\frac{1}{2\beta}}_t L^p_{x}}
\lesssim&\ \left\|B_{f^\omega} \right\|_{L^{\frac{1}{2\beta}}_tL^p_x}^2.
\end{split}
\end{equation}

\medskip

Obviously $\Phi_2$ and $\Phi_3$ can be estimated similarly. Moreover, we observe that $\|\Phi_2\|_{\mathcal H^{\frac{7}{2}-2\alpha}(\mathbb T^3)}$ shares an analogous estimate with $\|\Phi_4\|_{\mathcal H^{\frac{7}{2}-2\alpha}(\mathbb T^3)}$, while $\|\Phi_2\|_{L^{(\beta, \frac{1}{\beta})}_t L^p_{x}}$ and $\|\Phi_2\|_{L^{\frac{1}{2\beta}}_t L^p_{x}}$ enjoy analogous estimates with $\|\Phi_1\|_{L^{(\beta, \frac{1}{\beta})}_t L^p_{x}}$ and $\|\Phi_1\|_{L^{\frac{1}{2\beta}}_t L^p_{x}}$ respectively. Therefore we claim 
\begin{equation}\label{est-phi6}
\begin{split}
&\|\Phi_2\|_{\mathcal H^{\frac{7}{2}-2\alpha}(\mathbb T^3)}+\|\Phi_3\|_{\mathcal H^{\frac{7}{2}-2\alpha}(\mathbb T^3)}\\
\lesssim&\ \|B_{f^\omega}\|_{L^{(\beta,\frac{1}{\beta})}_tL^p_x} \|H\|_{L^{(\beta,\frac{1}{\beta})}_tL^p_x}
+\|B_{f^\omega}\|_{L^{(\beta,\frac{1}{\beta})}_tL^4_x}\|H\|_{L^{(\beta,\frac{1}{\beta})}_tL^p_x}\\
&+\|B_{f^\omega}\|_{L^{(\eta, 2)}_t \mathcal H^{\frac{11}{2}-2\alpha}_x} \|H\|_{L^{(\beta, \frac{1}{\beta})}_t L^{p}_x},\\
&\|\Phi_2\|_{L^{(\beta, \frac{1}{\beta})}_t L^p_{x}}+\|\Phi_3\|_{L^{(\beta, \frac{1}{\beta})}_t L^p_{x}}\\
\lesssim&\ \left\|B_{f^\omega} \right\|_{L^{\frac{1}{2\beta}}_tL^p_x}\left\|H \right\|_{L^{\frac{1}{2\beta}}_tL^p_x}+\left\|B_{f^\omega} \right\|_{L^{\frac{1}{2\beta}}_tL^p_x}\left\|H \right\|_{L^{(\beta, \frac{1}{\beta})}_t L^p_{x}},\\
&\|\Phi_2\|_{L^{\frac{1}{2\beta}}_t L^p_{x}}+\|\Phi_3\|_{L^{\frac{1}{2\beta}}_t L^p_{x}}
\lesssim  \left\|B_{f^\omega}\right\|_{L^{\frac{1}{2\beta}}_tL^p_x}\left\|H \right\|_{L^{\frac{1}{2\beta}}_tL^p_x}.
\end{split}
\end{equation}
It follows from the estimates (\ref{est-phi1})-(\ref{est-phi6}) and Lemma \ref{le-heat7} that for almost every $\omega\in \Omega$
\begin{equation}\label{map1}
\|\Phi(H)\|_{\mathcal Y}\lesssim \|H\|_{\mathcal Y}^2+\lambda^2
\end{equation}
which verifies condition (i) that $\Phi$ maps the subspace $\mathcal Y$ onto $\mathcal Y$. 

Regarding (ii), straightforward algebra shows that for $H_1, H_2\in\mathcal Y$
\begin{equation}\notag
\begin{split}
\Phi(H_1)(t)-\Phi(H_2)(t)
=&-\int_0^t e^{-(t-\tau)(-\Delta)^\alpha} \nabla\times\nabla\cdot (H_1\otimes (H_1-H_2))\, d\tau\\
&-\int_0^t e^{-(t-\tau)(-\Delta)^\alpha} \nabla\times\nabla\cdot ((H_1-H_2)\otimes H_2)\, d\tau\\
&-\int_0^t e^{-(t-\tau)(-\Delta)^\alpha} \nabla\times\nabla\cdot (B_{f^\omega}\otimes (H_1-H_2))\, d\tau\\
&-\int_0^t e^{-(t-\tau)(-\Delta)^\alpha} \nabla\times\nabla\cdot ((H_1-H_2)\otimes B_{f^\omega})\, d\tau.
\end{split}
\end{equation}
One notices that $\Phi(H_1)-\Phi(H_2)$ can be estimated similarly as $\Phi(H)$ in the space $\mathcal Y$ and 
\begin{equation}\label{map2}
\|\Phi(H_1)-\Phi(H_2)\|_{\mathcal Y}\lesssim \left(\|H_1\|_{\mathcal Y}+\|H_2\|_{\mathcal Y}+\lambda\right)\|H_1-H_2\|_{\mathcal Y}.
\end{equation}

We are ready to finish the proof the theorem by applying a fixed point argument. Indeed, following from (\ref{map1}) and (\ref{map2}), there exists a constant $C>0$ such that 
\begin{equation}\label{map3}
\begin{split}
\|\Phi(H)\|_{\mathcal Y}\leq&\ C\left( \|H\|_{\mathcal Y}^2+\lambda^2\right),\\
\|\Phi(H_1)-\Phi(H_2)\|_{\mathcal Y}\leq&\ C \left(\|H_1\|_{\mathcal Y}+\|H_2\|_{\mathcal Y}+\lambda\right)\|H_1-H_2\|_{\mathcal Y}.
\end{split}
\end{equation}
For such $C$, we then choose $\lambda$ such that 
\begin{equation}\notag
C\left((2C\lambda^2)^2+\lambda^2\right)\leq 2C\lambda^2, \ \ C\left(4C\lambda^2+\lambda\right)<1
\end{equation}
which are satisfied for $C^2\lambda^2\leq \frac18$.  Thus a suitable choice is $\lambda=\bar\lambda=\frac{1}{3 C}$. Therefore thanks to (\ref{map3}), for such $C$ and $\bar\lambda$, the map $\Phi$ is a contraction on the ball $B(0, 2C\bar\lambda^2)=B(0, \frac{2}{9C})\subset \mathcal Y$. Denote 
\begin{equation}\notag
\begin{split}
&E(f,\beta,\alpha, \bar\lambda, T)\\
=&\ \left\{\omega\in\Omega: \|B_{f^\omega}\|_{L^{\frac{1}{2\beta}}_TL^p_x}+\|B_{f^\omega}\|_{L^{(\beta, \frac{1}{\beta})}_TL^p_x}+\|B_{f^\omega}\|_{L^{(\eta,2)}_T\mathcal H^{\frac{11}{2}-2\alpha}_x}\geq \bar\lambda\right\}.
\end{split}
\end{equation}
Take $\Omega_T=E^c(f,\beta,\alpha, \bar\lambda, T)$ and $\Sigma=\cup_{-\infty<j<\infty}\Omega_{2^j}$. One can see that 
\[E(f,\beta,\alpha, \bar\lambda, T_1)\subset E(f,\beta,\alpha, \bar\lambda, T_2)\ \ \ \mbox{for}\ \ \ T_1\leq T_2\]
and hence $\Omega_{T_2}\subset \Omega_{T_1}$. For any $T>0$ there exists $k$ such that $2^{k-1}\leq T\leq 2^k$ and
\[\Omega_{2^k}\subset \Omega_{T}\subset \Omega_{2^{k-1}}.\]
Therefore we infer from Lemma \ref{le-heat7} that
\begin{equation}\notag
P(\Omega_T)\geq P(\Omega_{2^k})=1-P(E(f,\beta,\alpha, \bar\lambda, 2^k))\geq 1- 3c_1e^{-\frac{c_2\bar\lambda^2}{\|f\|^2_{\mathcal H^s_x}}}. 
\end{equation}
It follows that for any $\delta\in(0,1)$, if 
\begin{equation}\label{small}
\|f\|^2_{\mathcal H^s_x}\leq \frac{c_2}{9C^2(\ln(3c_1)-\ln(1-\delta))},
\end{equation} 
we have $P(\Omega_T)\geq \delta$ and hence $P(\Sigma)\geq \delta$.
On the other hand, without any assumption on the size of $\|f\|_{\mathcal H^s_x}$, (\ref{small}) indicates
$P(\Omega_T)\geq \delta$ holds for 
\begin{equation}\label{small2}
C^2\leq \frac{c_2}{9\|f\|^2_{\mathcal H^s_x} (\ln(3c_1)-\ln(1-\delta))}. 
\end{equation}
In view of the integral form (\ref{H2}) and $H(x,0)=0$, we know the constant $C$ appeared in (\ref{map3}) satisfies (\ref{small2}) for small time $T>0$. Therefore, we conclude the almost sure local well-posedness of (\ref{H}) in the space $\mathcal Y$ for general initial data $f\in \mathcal H^s$, and almost sure global well-posedness for small $f\in \mathcal H^s$.

\cbdu

\bigskip

\section{Well-posedness of the hyperdissipative NSE}\label{sec-nse}

This section concerns the well-posedness for the hyperdissipative NSE (\ref{hnse}). Namely we will prove Theorem \ref{thm-nse}. The proof is analogous to that of Theorem \ref{thm1} presented in Section \ref{sec-emhd}. We only include limited details to reveal the requirement on the parameters as stated in Theorem \ref{thm-nse}.

We consider a solution to (\ref{hnse}) in the form 
$u=u_{g^\omega}+V$
with $u_{g^\omega}=e^{-t(-\Delta)^\alpha} g^\omega$ and $V$ satisfying 
\begin{equation}\label{V}
\begin{split}
V_t+\left(u_{g^\omega}+V\right)\cdot\nabla \left(u_{g^\omega}+V\right)+\nabla\Pi=&-(-\Delta)^\alpha V,\\
\nabla\cdot V=&\ 0,\\
V(x,0)=&\ 0.
\end{split}
\end{equation}
Denote the map through the integral form of (\ref{V})
\begin{equation}\notag
\Psi(V)(t)=-\int_0^t e^{-(t-\tau)^{2\alpha}} \mathbb P \nabla\cdot \left(\left(u_{g^\omega}+V\right)\otimes \left(u_{g^\omega}+V\right)\right)\, d\tau.
\end{equation}
Define the subspace $\mathcal X\subset C([0,T]; \mathcal H^{\frac52-2\alpha}(\mathbb T^3))$ as
\begin{equation}\notag
\mathcal X= C([0,T]; \mathcal H^{\frac52-2\alpha}(\mathbb T^3))\cap L^{(\gamma, \frac1\gamma)}(0, T; L^q(\mathbb T^3))\cap L^{\frac1{2\gamma}}(0, T; L^q(\mathbb T^3))
\end{equation}
where $\gamma>0$ and $q\geq 2$ satisfy $\frac1{2\alpha}+\frac{3}{2q\alpha}+2\gamma=1$ as in Lemma \ref{le-heat9}.
We only need to show that 
\begin{equation}\label{psi}
\begin{split}
\|\Psi(V)\|_{\mathcal X}\lesssim&\ \lambda^2+\|V\|_{\mathcal X}^2 \ \ \ \forall \ \ V\in \mathcal X,\\
\|\Psi(V_1)-\Psi(V_2)\|_{\mathcal X}\lesssim&\ \left(\lambda+\|V_1\|_{\mathcal X}+\|V_2\|_{\mathcal X}\right)\|V_1-V_2\|_{\mathcal X} \ \ \ \forall \ \ V_1, V_2\in \mathcal X.
\end{split}
\end{equation}
Then a similar probability analysis combined with the fixed point argument as in the previous section provides a proof of Theorem \ref{thm-nse}.

Note that it is sufficient to prove the first inequality of (\ref{psi}) since the second one can be obtained similarly. Moreover, we only show details to estimate terms involving $V\otimes V$ and $u_{g^\omega}\otimes u_{g^\omega}$ as the mixed terms $u_{g^\omega}\otimes V$ and $V\otimes u_{g^\omega}$ can be estimated analogously. 

We first estimate $\int_0^t e^{-(t-\tau)(-\Delta)^\alpha} \mathbb P\nabla\cdot (V(x, \tau)\otimes V(x,\tau))\, d\tau$ in $\mathcal H^{\frac{5}{2}-2\alpha}(\mathbb T^3)$,
\begin{equation}\notag
\begin{split}
&\left\|\int_0^t e^{-(t-\tau)(-\Delta)^\alpha} \mathbb P\nabla\cdot (V(x, \tau)\otimes V(x,\tau))\, d\tau\right\|_{\mathcal H^{\frac{5}{2}-2\alpha}(\mathbb T^3)}\\
=&\ \left\| (1-\Delta)^{\frac54-\alpha}\int_0^t e^{-(t-\tau)(-\Delta)^\alpha} \mathbb P\nabla\cdot (V(x, \tau)\otimes V(x,\tau))\, d\tau \right\|_{L^2_x}\\
%=&\ \left\| \int_0^t e^{-(t-\tau)(-\Delta)^\alpha}(-\Delta)^{\frac12} (1-\Delta)^{\frac54-\alpha}\frac{\mathbb P\nabla\cdot}{(-\Delta)^{\frac12}} (V(x, \tau)\otimes V(x,\tau))\, d\tau \right\|_{L^2_x}\\
\lesssim&\  \int_0^t \left\|e^{-(t-\tau)(-\Delta)^\alpha}(-\Delta)^{\frac12} (1-\Delta)^{\frac54-\alpha} (V(x, \tau)\otimes V(x,\tau)) \right\|_{L^2_x}\, d\tau \\
\lesssim&\  \int_0^t \left\|e^{-(t-\tau)|\xi|^{2\alpha}}|\xi| \mathcal F\left((1-\Delta)^{\frac54-\alpha} (V(x, \tau)\otimes V(x,\tau))\right) \right\|_{L^2_{\xi}}\, d\tau 
\end{split}
\end{equation}
%where we used Plancherel's theorem in the last step. 
followed by applying H\"older's inequality and Lemma \ref{le-heat1} 
\begin{equation}\notag
\begin{split}
&\left\|e^{-(t-\tau)|\xi|^{2\alpha}}|\xi| \mathcal F\left((1-\Delta)^{\frac54-\alpha} (V(x, \tau)\otimes V(x,\tau))\right) \right\|_{L^2_{\xi}}\\
\lesssim&\ \left\|e^{-(t-\tau)|\xi|^{2\alpha}}|\xi|\right\|_{L^q_{\xi}} \left\|\mathcal F\left((1-\Delta)^{\frac54-\alpha} (V(x, \tau)\otimes V(x,\tau))\right) \right\|_{L^{\frac{2q}{q-2}}_{\xi}}\\
\lesssim&\ (t-\tau)^{-\frac{1}{2\alpha}-\frac{3}{2q\alpha}} \left\|(1-\Delta)^{\frac54-\alpha} (V(x, \tau)\otimes V(x,\tau)) \right\|_{L^{\frac{2q}{q+2}}_{x}}\\
\lesssim&\ (t-\tau)^{-\frac{1}{2\alpha}-\frac{3}{2q\alpha}} \left\| V\nabla^{\frac52-2\alpha}V \right\|_{L^{\frac{2q}{q+2}}_{x}}\\
\lesssim&\ (t-\tau)^{-\frac{1}{2\alpha}-\frac{3}{2q\alpha}} \left\| V \right\|_{L^{q}_{x}} \left\| \nabla^{\frac52-2\alpha}V \right\|_{L^{2}_{x}}.
\end{split}
\end{equation}
Combining the last two inequalities yields
\begin{equation}\notag
\begin{split}
&\|\Phi_1\|_{\mathcal H^{\frac{5}{2}-2\alpha}(\mathbb T^3)}\\
\lesssim&\ \int_0^t (t-\tau)^{-\frac{1}{2\alpha}-\frac{3}{2q\alpha}} \left\| V \right\|_{L^{q}_{x}} \left\| \nabla^{\frac52-2\alpha}V \right\|_{L^{2}_{x}}\, d\tau\\
\lesssim&\ \left\|V \right\|_{L^\infty_t\mathcal H^{\frac52-2\alpha}_{x}}\int_0^t (t-\tau)^{-\frac{1}{2\alpha}-\frac{3}{2q\alpha}}\tau^{-\gamma} \left(\tau^\gamma\left\| V \right\|_{L^{q}_{x}}\right) \, d\tau\\
\lesssim&\ \left\|V \right\|_{L^\infty_t\mathcal H^{\frac52-2\alpha}_{x}}\left(\int_0^t (t-\tau)^{\left(-\frac{1}{2\alpha}-\frac{3}{2q\alpha}\right)\frac{1}{1-\gamma}}\tau^{-\frac{\gamma}{1-\gamma}} \, d\tau\right)^{1-\gamma} \left(\int_0^t \tau \| V \|_{L^{q}_{x}}^{\frac1\gamma}\, d\tau\right)^{\gamma}.
\end{split}
\end{equation}
Noticing that 
\[\left(\frac{1}{2\alpha}+\frac{3}{2q\alpha}\right)\frac{1}{1-\gamma}+\frac{\gamma}{1-\gamma}=1,\]
\[0<\left(\frac{1}{2\alpha}+\frac{3}{2q\alpha}\right)\frac{1}{1-\gamma}<1, \ \ 0<\frac{\gamma}{1-\gamma}<1,\]
It follows from Lemma \ref{le-beta} that
\begin{equation}\notag
\int_0^t (t-\tau)^{\left(-\frac{1}{2\alpha}-\frac{3}{2q\alpha}\right)\frac{1}{1-\gamma}}\tau^{-\frac{\gamma}{1-\gamma}} \, d\tau
\lesssim 1.
\end{equation}
Therefore we have
\begin{equation}\label{est-psi1}
\begin{split}
&\left\|\int_0^t e^{-(t-\tau)(-\Delta)^\alpha} \mathbb P\nabla\cdot (V(x, \tau)\otimes V(x,\tau))\, d\tau\right\|_{\mathcal H^{\frac{5}{2}-2\alpha}(\mathbb T^3)}\\
\lesssim&\ \left\|V \right\|_{L^\infty_t\mathcal H^{\frac52-2\alpha}_{x}}\left\|V \right\|_{L^{(\gamma, \frac{1}{\gamma})}_t L^q_{x}}.
\end{split}
\end{equation}

The estimate of $\| \int_0^t e^{-(t-\tau)(-\Delta)^\alpha} \mathbb P\nabla\cdot (V(x, \tau)\otimes V(x,\tau))\, d\tau \|_{L^{(\gamma, \frac{1}{\gamma})}_t L^q_{x}}$ is given by
\begin{equation}\notag
\begin{split}
&\left\|t^\gamma \int_0^t e^{-(t-\tau)(-\Delta)^\alpha} \mathbb P\nabla\cdot (V(x, \tau)\otimes V(x,\tau))\, d\tau \right\|_{L^{\frac{1}{\gamma}}_t L^q_{x}}\\
\lesssim &\ \left\|t^\gamma \int_0^t \left\|e^{-(t-\tau)(-\Delta)^\alpha} (-\Delta)^{\frac12} (V(x, \tau)\otimes V(x,\tau))\right\|_{L^q_{x}}\, d\tau \right\|_{L^{\frac{1}{\gamma}}_t}
\end{split}
\end{equation}
and %Applying H\"older's inequality and Lemma \ref{le-heat1} yields
\begin{equation}\notag
\begin{split}
&\left\|e^{-(t-\tau)(-\Delta)^\alpha} (-\Delta)^{\frac12} (V(x, \tau)\otimes V(x,\tau))\right\|_{L^q_{x}}\\
\lesssim&\ \left\|e^{-(t-\tau)|\xi|^{2\alpha}} |\xi| \mathcal F\left((V(x, \tau)\otimes V(x,\tau))\right)\right\|_{L^{\frac{q}{q-1}}_{\xi}}\\
\lesssim&\ \left\|e^{-(t-\tau)|\xi|^{2\alpha}} |\xi|\right\|_{L^{q}_{\xi}}\left\| \mathcal F\left((V(x, \tau)\otimes V(x,\tau))\right)\right\|_{L^{\frac{q}{q-2}}_{\xi}}\\
\lesssim&\ (t-\tau)^{-\frac1{2\alpha}-\frac{3}{2q\alpha}}\left\|V(x, \tau)\otimes V(x,\tau)\right\|_{L^{\frac{q}{2}}_{x}}\\
\lesssim&\ (t-\tau)^{-\frac1{2\alpha}-\frac{3}{2q\alpha}}\left\|V\right\|^2_{L^{q}_{x}}.
\end{split}
\end{equation}
Therefore it follows
\begin{equation}\notag
\begin{split}
&\left\| \int_0^t e^{-(t-\tau)(-\Delta)^\alpha} \mathbb P\nabla\cdot (V(x, \tau)\otimes V(x,\tau))\, d\tau \right\|_{L^{(\gamma, \frac{1}{\gamma})}_t L^q_{x}}\\
\lesssim &\ \left\|\int_0^t t^\gamma(t-\tau)^{-\frac1{2\alpha}-\frac{3}{2q\alpha}}\left\|V\right\|^2_{L^{q}_{x}} \, d\tau \right\|_{L^{\frac{1}{\gamma}}_t}\\
\lesssim &\ \left\| \int_0^{\frac{t}{2}} t^\gamma (t-\tau)^{-\frac1{2\alpha}-\frac{3}{2q\alpha}}\left\|V\right\|^2_{L^{q}_{x}} \, d\tau \right\|_{L^{\frac{1}{\gamma}}_t}\\
&+ \left\|\int_{\frac{t}{2}}^t t^\gamma(t-\tau)^{-\frac1{2\alpha}-\frac{3}{2q\alpha}}\left\|V\right\|^2_{L^{q}_{x}} \, d\tau \right\|_{L^{\frac{1}{\beta}}_t}\\
=:&\ I_5+I_6.
\end{split}
\end{equation}
The term $I_5$ is estimated as by using H\"older's inequality
\begin{equation}\notag
\begin{split}
& \left\| \int_0^{\frac{t}{2}} t^\gamma (t-\tau)^{-\frac1{2\alpha}-\frac{3}{2q\alpha}}\left\|V\right\|^2_{L^{q}_{x}} \, d\tau \right\|_{L^{\frac{1}{\gamma}}_t}\\
\lesssim&\ \left\| \int_0^{\frac{t}{2}}  (t-\tau)^{\gamma-\frac1{2\alpha}-\frac{3}{2q\alpha}}\left\|V\right\|^2_{L^{q}_{x}} \, d\tau \right\|_{L^{\frac{1}{\gamma}}_t}\\
\lesssim&\ \left(\int_0^{\frac{t}{2}} (t-\tau)^{(\gamma-\frac{1}{2\alpha}-\frac{3}{2q\alpha})(\frac{1}{1-3\gamma})}\, d\tau\right)^{1-3\gamma}
\left(\int_0^{\frac{t}{2}}\left\|V\right\|^{2\cdot \frac{1}{4\gamma}}_{L^{q}_{x}}\, d\tau \right)^{4\gamma}\\
=&\ \left(\int_0^{\frac{t}{2}} (t-\tau)^{-1}\, d\tau\right)^{1-3\gamma} \|V\|^2_{L^{\frac{1}{2\gamma}}_t L^{q}_{x}}\\
\lesssim&\ \|V\|^2_{L^{\frac{1}{2\gamma}}_t L^{q}_{x}}
\end{split}
\end{equation}
and the term $I_6$ is estimated by applying H\"older's inequality and Hardy-Littlewood-Sobolev lemma
\begin{equation}\notag
\begin{split}
& \left\| \int_{\frac{t}{2}}^t t^\gamma (t-\tau)^{-\frac1{2\alpha}-\frac{3}{2q\alpha}}\left\|H\right\|^2_{L^{q}_{x}} \, d\tau \right\|_{L^{\frac{1}{\gamma}}_t}\\
\lesssim&\ \left\| \int_0^{\frac{t}{2}}  (t-\tau)^{-\frac1{2\alpha}-\frac{3}{2q\alpha}}\tau^\gamma \|V\|_{L^{q}_{x}} \|V\|_{L^{q}_{x}}\, d\tau \right\|_{L^{\frac{1}{\gamma}}_t}\\
\lesssim&\ \left\| \left(\int_{\frac{t}{2}}^t\tau \|V\|_{L^{q}_{x}}^{\frac{1}{\gamma}} \, d\tau\right)^{\gamma}\left(\int_{\frac{t}{2}}^t  (t-\tau)^{(-\frac1{2\alpha}-\frac{3}{2q\alpha})\frac{1}{1-\gamma}} \|V\|_{L^{q}_{x}}^{\frac{1}{1-\gamma}}\, d\tau\right)^{1-\gamma} \right\|_{L^{\frac{1}{\gamma}}_t}\\
\lesssim&\ \|V\|_{L^{(\gamma, \frac{1}{\gamma})}_t L^{q}_{x}}\left\| \int_{\frac{t}{2}}^t  (t-\tau)^{(-\frac1{2\alpha}-\frac{3}{2q\alpha})\frac{1}{1-\gamma}} \|V\|_{L^{q}_{x}}^{\frac{1}{1-\gamma}}\, d\tau \right\|^{1-\gamma}_{L^{\frac{1-\gamma}{\gamma}}_t}\\
\lesssim&\  \|V\|_{L^{(\gamma, \frac{1}{\gamma})}_t L^{q}_{x}} 
\left\| \|V\|_{L^{q}_{x}}^{\frac{1}{1-\gamma}} \right\|_{L^{\frac{1-\gamma}{2\gamma}}_t}^{1-\gamma}\\
\lesssim&\ \|V\|_{L^{(\gamma, \frac{1}{\gamma})}_t L^{q}_{x}}\|V\|_{L^{\frac{1}{2\gamma}}_t L^{q}_{x}}.
\end{split}
\end{equation}
The estimates above together imply
\begin{equation}\label{est-psi2}
\begin{split}
&\|\int_0^t e^{-(t-\tau)(-\Delta)^\alpha} \mathbb P\nabla\cdot (V(x, \tau)\otimes V(x,\tau))\, d\tau\|_{L^{(\gamma, \frac{1}{\gamma})}_t L^q_{x}}\\
\lesssim&\ \left\|V \right\|_{L^{\frac{1}{2\gamma}}_tL^q_x}^2+\left\|V \right\|_{L^{(\gamma, \frac{1}{\gamma})}_t L^q_{x}}\left\|V \right\|_{L^{\frac{1}{2\gamma}}_tL^q_x}.
\end{split}
\end{equation}

Finally we estimate $\| \int_0^t e^{-(t-\tau)(-\Delta)^\alpha} \mathbb P\nabla\cdot (V(x, \tau)\otimes V(x,\tau))\, d\tau\|_{L^{\frac{1}{2\gamma}}_t L^q_{x}}$, 
\begin{equation}\notag
\begin{split}
&\left\| \int_0^t e^{-(t-\tau)(-\Delta)^\alpha} \mathbb P\nabla\cdot (V(x, \tau)\otimes V(x,\tau))\, d\tau \right\|_{L^{\frac{1}{2\gamma}}_t L^q_{x}}\\
\lesssim &\ \left\| \int_0^t \left\|e^{-(t-\tau)(-\Delta)^\alpha} (-\Delta)^{\frac12} (V(x, \tau)\otimes V(x,\tau))\right\|_{L^q_{x}}\, d\tau \right\|_{L^{\frac{1}{2\gamma}}_t},
\end{split}
\end{equation}
and since %Similarly as before, we deduce by applying H\"older's inequality and Lemma \ref{le-heat1} 
\begin{equation}\notag
\begin{split}
& \left\|e^{-(t-\tau)(-\Delta)^\alpha} (-\Delta)^{\frac12} (V(x, \tau)\otimes V(x,\tau))\right\|_{L^q_{x}}\\
\lesssim &\ \left\|e^{-(t-\tau)|\xi|^{2\alpha}} |\xi| \mathcal F\left(V(x, \tau)\otimes V(x,\tau)\right)\right\|_{L^{\frac{q}{q-1}}_{\xi}}\\
\lesssim &\ \left\|e^{-(t-\tau)|\xi|^{2\alpha}} |\xi|\right\|_{L^{q}_{\xi}} \left\|\mathcal F\left(V(x, \tau)\otimes V(x,\tau)\right)\right\|_{L^{\frac{q}{q-2}}_{\xi}}\\
\lesssim &\ (t-\tau)^{-\frac{1}{2\alpha}-\frac{3}{2q\alpha}} \left\|V(x, \tau)\otimes V(x,\tau)\right\|_{L^{\frac{q}{2}}_{x}}\\
\lesssim &\ (t-\tau)^{-\frac{1}{2\alpha}-\frac{3}{2q\alpha}} \|V\|^2_{L^{q}_{x}}.
\end{split}
\end{equation}
it follows from Hardy-Littlewood-Sobolev lemma that
\begin{equation}\label{est-psi3}
\begin{split}
&\left\|\int_0^t e^{-(t-\tau)(-\Delta)^\alpha} \mathbb P\nabla\cdot (V(x, \tau)\otimes V(x,\tau))\, d\tau\right\|_{L^{\frac{1}{2\gamma}}_t L^q_{x}}\\
\lesssim&\ \left\| \int_0^t (t-\tau)^{-\frac{1}{2\alpha}-\frac{3}{2q\alpha}} \|V\|^2_{L^{q}_{x}} \, d\tau \right\|_{L^{\frac{1}{2\gamma}}_t}\\
%\lesssim&\ \left\| \|V\|^2_{L^{q}_{x}}\right\|_{L^{\frac{1}{4\gamma}}_t}\\
\lesssim&\  \left\|V \right\|_{L^{\frac{1}{2\gamma}}_tL^q_x}^2.
\end{split}
\end{equation}
Therefore it follows from (\ref{est-psi1})-(\ref{est-psi3}) that 
\begin{equation}\label{est-psi4}
\left\|\int_0^t e^{-(t-\tau)(-\Delta)^\alpha} \mathbb P\nabla\cdot (V(x, \tau)\otimes V(x,\tau))\, d\tau\right\|_{\mathcal X}
\lesssim  \left\|V \right\|_{\mathcal X}^2.
\end{equation}

Now we estimate  $\int_0^t e^{-(t-\tau)(-\Delta)^\alpha} \mathbb P\nabla\cdot (u_{g^\omega}\otimes u_{g^\omega})\, d\tau$ in $\mathcal X$. First, we have
\begin{equation}\notag
\begin{split}
&\left\|\int_0^t e^{-(t-\tau)(-\Delta)^\alpha} \mathbb P\nabla\cdot (u_{g^\omega}\otimes u_{g^\omega})\, d\tau\right\|_{\mathcal H^{\frac{5}{2}-2\alpha}(\mathbb T^3)}\\
=&\ \left\| (1-\Delta)^{\frac54-\alpha}\int_0^t e^{-(t-\tau)(-\Delta)^\alpha} \mathbb P\nabla\cdot (u_{g^\omega}(x, \tau)\otimes u_{g^\omega}(x,\tau))\, d\tau \right\|_{L^2_x}\\
\lesssim&\  \int_0^\delta \left\|e^{-(t-\tau)|\xi|^{2\alpha}}|\xi|^{1+\frac52-2\alpha} \mathcal F\left(u_{g^\omega}(x, \tau)\otimes u_{g^\omega}(x,\tau)\right) \right\|_{L^2_{\xi}}\, d\tau \\
&+ \int_\delta^t \left\|e^{-(t-\tau)|\xi|^{2\alpha}} \mathcal F\left((-\Delta)^{\frac12}(1-\Delta)^{\frac54-\alpha} (u_{g^\omega}(x, \tau)\otimes u_{g^\omega}(x,\tau))\right) \right\|_{L^2_{\xi}}\, d\tau \\
=:&\ I_7+I_8
\end{split}
\end{equation}
for a small constant $0<\delta<t$. For $2\leq q\leq 4$, we have from H\"older's inequality and Lemma \ref{le-heat1} 
\begin{equation}\notag
\begin{split}
&\left\|e^{-(t-\tau)|\xi|^{2\alpha}}|\xi|^{1+\frac52-2\alpha} \mathcal F\left(u_{g^\omega}(x, \tau)\otimes u_{g^\omega}(x,\tau)\right) \right\|_{L^2_{\xi}}\\
\lesssim &\ \left\|e^{-(t-\tau)|\xi|^{2\alpha}}|\xi|^{1+\frac52-2\alpha}\right\|_{L^{\frac{2q}{4-q}}_{\xi}}\left\|\mathcal F\left(u_{g^\omega}(x, \tau)\otimes u_{g^\omega}(x,\tau)\right) \right\|_{L^{\frac{q}{q-2}}_{\xi}}\\
\lesssim &\ (t-\tau)^{-\frac{7}{4\alpha}+1-\frac{3}{2\alpha}\cdot\frac{4-q}{2q}}\left\|u_{g^\omega}(x, \tau)\otimes u_{g^\omega}(x,\tau) \right\|_{L^{\frac{q}{2}}_{x}}\\
\lesssim &\ (t-\tau)^{-\frac{7}{4\alpha}+1-\frac{3(4-q)}{4q\alpha}}\|u_{g^\omega} \|^2_{L^{q}_{x}},
\end{split}
\end{equation}
and hence
\begin{equation}\notag
\begin{split}
I_7\lesssim&\ \int_0^\delta (t-\tau)^{-\frac{7}{4\alpha}+1-\frac{3(4-q)}{4q\alpha}}\|u_{g^\omega} \|^2_{L^{q}_{x}}\, d\tau\\
=&\ \int_0^\delta \left((t-\tau)^{-\frac{7}{4\alpha}+1-\frac{3(4-q)}{4q\alpha}}\tau^{-2\gamma}\right) \left(\tau^{\gamma}\|u_{g^\omega} \|_{L^{q}_{x}}\right)^2\, d\tau\\
\lesssim&\ \left(\int_0^\delta (t-\tau)^{(-\frac{7}{4\alpha}+1-\frac{3(4-q)}{4q\alpha})\frac{1}{1-2\gamma}}\tau^{-2\gamma\cdot\frac{1}{1-2\gamma}} \, d\tau\right)^{1-2\gamma} \|u_{g^\omega}\|^2_{L^{(\gamma,\frac{1}{\gamma})}_tL^q_x}\\
=&\ \left(\int_0^\delta (t-\tau)^{-1+\frac{2\gamma}{1-2\gamma}}\tau^{-\frac{2\gamma}{1-2\gamma}} \, d\tau\right)^{1-2\gamma} \|u_{g^\omega}\|^2_{L^{(\gamma,\frac{1}{\gamma})}_tL^q_x}\\
\lesssim&\ \|u_{g^\omega}\|^2_{L^{(\gamma,\frac{1}{\gamma})}_tL^q_x}.
\end{split}
\end{equation}
For $q>4$, we instead have
\begin{equation}\notag
\begin{split}
&\left\|e^{-(t-\tau)|\xi|^{2\alpha}}|\xi|^{1+\frac52-2\alpha} \mathcal F\left(u_{g^\omega}(x, \tau)\otimes u_{g^\omega}(x,\tau)\right) \right\|_{L^2_{\xi}}\\
\lesssim &\ \left\|e^{-(t-\tau)|\xi|^{2\alpha}}|\xi|^{1+\frac52-2\alpha}\right\|_{L^{\infty}_{\xi}}\left\|\mathcal F\left(u_{g^\omega}(x, \tau)\otimes u_{g^\omega}(x,\tau)\right) \right\|_{L^{2}_{\xi}}\\
\lesssim &\ (t-\tau)^{-\frac{7}{4\alpha}+1}\left\|u_{g^\omega}(x, \tau)\otimes u_{g^\omega}(x,\tau) \right\|_{L^{2}_{x}}\\
\lesssim &\ (t-\tau)^{-\frac{7}{4\alpha}+1}\|u_{g^\omega} \|^2_{L^{4}_{x}}
\end{split}
\end{equation}
and 
\begin{equation}\notag
\begin{split}
I_7\lesssim&\ \int_0^\delta (t-\tau)^{-\frac{7}{4\alpha}+1}\|u_{g^\omega} \|^2_{L^{4}_{x}}\, d\tau\\
=&\ \int_0^\delta \left((t-\tau)^{-\frac{7}{4\alpha}+1}\tau^{-2\gamma}\right) \left(\tau^{\gamma}\|u_{g^\omega} \|_{L^{4}_{x}}\right)^2\, d\tau\\
\lesssim&\ \left(\int_0^\delta (t-\tau)^{(-\frac{7}{4\alpha}+1)\frac{1}{1-2\gamma}}\tau^{-\frac{2\gamma}{1-2\gamma}} \, d\tau\right)^{1-2\gamma} \|u_{g^\omega}\|^2_{L^{(\gamma,\frac{1}{\gamma})}_tL^4_x}\\
\lesssim&\ \|u_{g^\omega}\|^2_{L^{(\gamma,\frac{1}{\gamma})}_tL^4_x}.
\end{split}
\end{equation}
%In the estimate above, the time integral is handled as
%\begin{equation}\notag
%\begin{split}
%&\int_0^\delta (t-\tau)^{(-\frac{7}{4\alpha}+1)\frac{1}{1-2\gamma}}\tau^{-\frac{2\gamma}{1-2\gamma}} \, d\tau\\
%< &\ t^{(-\frac{7}{4\alpha}+1)\frac{1}{1-2\gamma}} \int_0^\delta \tau^{-\frac{2\gamma}{1-2\gamma}} \, d\tau\\
%\lesssim&\ 1
%\end{split}
%\end{equation}
%since $t>\delta>0$ and observing that for $1<\alpha<\frac74$
%\[ (-\frac{7}{4\alpha}+1)\frac{1}{1-2\gamma}<0,\ \ \  0<\frac{2\gamma}{1-2\gamma}<1.\]
The term $I_8$ is estimated as % again invoking H\"older's inequality and Lemma \ref{le-heat1} we obtain
\begin{equation}\notag
\begin{split}
&\left\|e^{-(t-\tau)|\xi|^{2\alpha}} \mathcal F\left((-\Delta)^{\frac12}(1-\Delta)^{\frac54-\alpha} (u_{g^\omega}(x, \tau)\otimes u_{g^\omega}(x,\tau))\right) \right\|_{L^2_{\xi}}\\
\lesssim&\ \left\|e^{-(t-\tau)|\xi|^{2\alpha}} \right\|_{L^{q}_{\xi}} \left\|\mathcal F\left((-\Delta)^{\frac12}(1-\Delta)^{\frac54-\alpha} (u_{g^\omega}(x, \tau)\otimes u_{g^\omega}(x,\tau))\right) \right\|_{L^{\frac{2q}{q-2}}_{\xi}}\\
\lesssim&\ (t-\tau)^{-\frac{3}{2q\alpha}} \left\|(-\Delta)^{\frac12}(1-\Delta)^{\frac54-\alpha} (u_{g^\omega}(x, \tau)\otimes u_{g^\omega}(x,\tau)) \right\|_{L^{\frac{2q}{q+2}}_{x}}\\
\lesssim&\ (t-\tau)^{-\frac{3}{2q\alpha}} \|\nabla^{\frac{7}{2}-2\alpha} u_{g^\omega}\|_{L^{2}_{x}}  \|u_{g^\omega}\|_{L^{q}_{x}}.
\end{split}
\end{equation}
 and %Therefore we have for $\frac1{\alpha}-\frac12<\zeta<\frac12-2\gamma$ and $1<\alpha<\frac74$
\begin{equation}\notag
\begin{split}
I_8\lesssim&\ \int_\delta^t (t-\tau)^{-\frac{3}{2q\alpha}} \|\nabla^{\frac{7}{2}-2\alpha} u_{g^\omega}\|_{L^{2}_{x}}  \|u_{g^\omega}\|_{L^{q}_{x}} \, d\tau\\
=&\ \int_\delta^t (t-\tau)^{-\frac{3}{2q\alpha}}\tau^{-\zeta-\gamma} \left(\tau^\zeta\|\nabla^{\frac{7}{2}-2\alpha} u_{g^\omega}\|_{L^{2}_{x}}\right)\left(\tau^{\gamma}  \|u_{g^\omega}\|_{L^{q}_{x}}\right) \, d\tau\\
\lesssim &\ \left(\int_\delta^t (t-\tau)^{-\frac{3}{2q\alpha}\cdot \frac{4}{3-4\gamma}}\tau^{(-\zeta-\gamma)\frac{4}{3-4\gamma}} \, d\tau\right)^{\frac{3-4\gamma}{4}}\|u_{g^\omega}\|_{L^{(\zeta, 4)}_t \mathcal H^{\frac{7}{2}-2\alpha}_x} \|u_{g^\omega}\|_{L^{(\gamma, \frac{1}{\gamma})}_t L^{q}_x}\\
\lesssim &\ \|u_{g^\omega}\|_{L^{(\zeta, 4)}_t \mathcal H^{\frac{7}{2}-2\alpha}_x} \|u_{g^\omega}\|_{L^{(\gamma, \frac{1}{\gamma})}_t L^{q}_x}.
\end{split}
\end{equation}
The time integral above is bounded for $\frac1{2\alpha}+\frac14\gamma\leq\zeta<1-\frac74\gamma$ and $1<\alpha<\frac74$. Indeed, we have
\begin{equation}\notag
\begin{split}
&\int_\delta^t (t-\tau)^{-\frac{3}{2q\alpha}\cdot \frac{4}{4-3\gamma}}\tau^{(-\zeta-\gamma)\frac{4}{4-3\gamma}} \, d\tau\\
=&\ t^{1- \frac{6}{q\alpha(4-3\gamma)}-\frac{4(\zeta+\gamma)}{4-3\gamma}}\int_{\delta/t}^{1} (1-\tau')^{- \frac{6}{q\alpha(4-3\gamma)}} (\tau')^{-\frac{4(\zeta+\gamma)}{4-3\gamma}} \, d\tau'\\
\lesssim&\ 1 \ \ \ \mbox{for} \ \ 0< \delta\leq t
\end{split}
\end{equation}
since 
\[1- \frac{6}{q\alpha(4-3\gamma)}-\frac{4(\zeta+\gamma)}{4-3\gamma}\leq 0, \ \ 0<\frac{6}{q\alpha(4-3\gamma)}<1, \ \ 0<\frac{4(\zeta+\gamma)}{4-3\gamma}<1.\]
%\[\frac{\frac12+\zeta-\frac1{\alpha}}{1-2\gamma}>0, \ \ 0<\frac{2(\zeta+\gamma)}{1-2\gamma}<1, \ \ 0<\frac{3}{q\alpha(1-2\gamma)}<1.\]
Combining the estimates above gives
\begin{equation}\label{est-psi4}
\begin{split}
&\left\|\int_0^t e^{-(t-\tau)(-\Delta)^\alpha} \mathbb P\nabla\cdot (u_{g^\omega}\otimes u_{g^\omega})\, d\tau\right\|_{\mathcal H^{\frac{5}{2}-2\alpha}(\mathbb T^3)}\\
\lesssim&\ \|u_{g^\omega}\|^2_{L^{(\gamma,\frac{1}{\gamma})}_tL^q_x}
+\|u_{g^\omega}\|^2_{L^{(\gamma,\frac{1}{\gamma})}_tL^4_x}\\
&+\|u_{g^\omega}\|_{L^{(\zeta, 4)}_t \mathcal H^{\frac{7}{2}-2\alpha}_x} \|u_{g^\omega}\|_{L^{(\gamma, \frac{1}{\gamma})}_t L^{q}_x}.
\end{split}
\end{equation}
The estimates of $\int_0^t e^{-(t-\tau)(-\Delta)^\alpha} \mathbb P\nabla\cdot (u_{g^\omega}\otimes u_{g^\omega})\, d\tau$ in $L^{(\gamma, \frac{1}{\gamma})}_t L^{q}_x$ and $L^{\frac{1}{2\gamma}}_t L^{q}_x$ 
\begin{equation}\label{est-psi5}
\begin{split}
&\left\|\int_0^t e^{-(t-\tau)(-\Delta)^\alpha} \mathbb P\nabla\cdot (u_{g^\omega}\otimes u_{g^\omega})\, d\tau\right\|_{L^{(\gamma, \frac{1}{\gamma})}_t L^q_{x}}\\
\lesssim&\ \left\|u_{g^\omega} \right\|_{L^{\frac{1}{2\gamma}}_tL^q_x}^2+\left\|u_{g^\omega} \right\|_{L^{(\gamma, \frac{1}{\gamma})}_t L^q_{x}}\left\|u_{g^\omega} \right\|_{L^{\frac{1}{2\gamma}}_tL^q_x},\\
&\left\|\int_0^t e^{-(t-\tau)(-\Delta)^\alpha} \mathbb P\nabla\cdot (u_{g^\omega}\otimes u_{g^\omega})\, d\tau\right\|_{L^{\frac{1}{2\gamma}}_t L^q_{x}}
\lesssim \left\|u_{g^\omega} \right\|_{L^{\frac{1}{2\gamma}}_tL^q_x}^2
\end{split}
\end{equation}
can be obtained in an analogy way with the estimate of $\Phi_1$ in the previous section.  Thus it follows from (\ref{est-psi4}), (\ref{est-psi5}) and Lemmas \ref{le-heat9}, \ref{le-heat11} and \ref{le-heat12} that
\begin{equation}\label{est-psi6}
\left\|\int_0^t e^{-(t-\tau)(-\Delta)^\alpha} \mathbb P\nabla\cdot (u_{g^\omega}\otimes u_{g^\omega})\, d\tau\right\|_{\mathcal X}\\
\lesssim \lambda^2.
\end{equation}
In view of (\ref{est-psi4}) and (\ref{est-psi6}), we conclude the proof of the first inequality in (\ref{psi}) is complete and so is the proof of Theorem \ref{thm-nse}.

\bigskip

\section{Well-posedness of the full system of Hall MHD}\label{sec-full}

%For given initial data $u(x,0)=g$ and $B(x,0)=f$, we denote their randomization by $g^\omega$ and $f^\omega$ respectively. We further denote the free evolution of $g^\omega$ and $f^\omega$ respectively by 
%\[u_{g^\omega}=e^{-t(-\Delta)^\alpha}g^\omega, \ \ B_{f^\omega}=e^{-t(-\Delta)^\alpha}f^\omega. \]
We address the well-posedness of the Hall MHD system (\ref{hmhda}) in this final section. 
As before, we seek a solution of (\ref{hmhda}) in the form
\[u=u_{g^\omega}+V, \ \ B=B_{f^\omega}+H \]
with $(V, H)$ satisfying the system 
\begin{equation}\label{VH}
\begin{split}
V_t+(u_{g^\omega}+V) \cdot\nabla (u_{g^\omega}+V) \\
-(B_{f^\omega}+H)\cdot\nabla (B_{f^\omega}+H)+\nabla \Pi=&-(-\Delta)^\alpha V,\\
H_t+(u_{g^\omega}+V)\cdot\nabla) (B_{f^\omega}+H)-(B_{f^\omega}+H)\cdot\nabla (u_{g^\omega}+V)\\
+\nabla\times\nabla\cdot\left((B_{f^\omega}+H)\otimes (B_{f^\omega}+H) \right)=&-(-\Delta)^\alpha H,\\
\nabla\cdot V= 0, \ \ \nabla\cdot H=&\ 0,\\
V(x,0)=H(x,0)=&\ 0.
\end{split}
\end{equation}
Without causing confusion, we use the same notations $\Phi$ and $\Psi$ to denote the maps
\begin{equation}\notag
\begin{split}
\Psi(V, H)(t)=&-\int_0^t e^{-(t-\tau)(-\Delta)^\alpha}\mathbb P\nabla\cdot\left[(u_{g^\omega}+V)\otimes (u_{g^\omega}+V)\right]\, d\tau\\
&+\int_0^t e^{-(t-\tau)(-\Delta)^\alpha}\mathbb P\nabla\cdot\left[(B_{f^\omega}+H)\otimes (B_{f^\omega}+H)\right]\, d\tau\\
:=&\ \Psi^1(V, H)+\Psi^2(V, H),\\
\Phi(V, H)(t)=&-\int_0^t e^{-(t-\tau)(-\Delta)^\alpha}\nabla\cdot\left[(u_{g^\omega}+V)\otimes (B_{f^\omega}+H)\right]\, d\tau\\
&+\int_0^t e^{-(t-\tau)(-\Delta)^\alpha}\nabla\cdot\left[(B_{f^\omega}+H)\otimes (u_{g^\omega}+V)\right]\, d\tau\\
&-\int_0^t e^{-(t-\tau)(-\Delta)^\alpha}\nabla\times\nabla\cdot\left[(B_{f^\omega}+H)\otimes (B_{f^\omega}+H)\right]\, d\tau\\
:=&\ \Phi^1(V, H)+\Phi^2(V, H)+\Phi^3(V, H).
\end{split}
\end{equation}
Define the functional spaces $\mathcal X$ and $\mathcal Y$ as in previous sections
\begin{equation}\notag
\begin{split}
\mathcal X=&\ C([0,T]; \mathcal H^{\frac52-2\alpha}(\mathbb T^3))\cap L^{(\gamma,\frac{1}{\gamma})}_tL^{q}_x\cap L^{\frac{1}{2\gamma}}_tL^{q}_x,\\
\mathcal Y=&\ C([0,T]; \mathcal H^{\frac72-2\alpha}(\mathbb T^3))\cap L^{(\beta,\frac{1}{\beta})}_tL^{p}_x\cap L^{\frac{1}{2\beta}}_tL^{p}_x
\end{split}
\end{equation}
with 
\begin{equation}\label{para}
\frac{1}{\alpha}+\frac{3}{2p\alpha}+2\beta=1, \ \ \frac{1}{2\alpha}+\frac{3}{2q\alpha}+2\gamma=1
\end{equation}
for $\beta, \gamma>0$ and $p, q\geq 2$.
%{\color{red} Note that if $\beta=\gamma$, it then follows $1+\frac3p=\frac3q$. If $2\leq p\leq \infty$, then $\frac{6}{5}\leq q\leq 3$. Is that good enough for the NSE?}

The existence and uniqueness of solution to (\ref{hmhda}) is a consequence of that the map $(\Psi, \Phi)$ is a contraction on $\mathcal X \times\mathcal Y$. Thus in order to prove Theorem \ref{thm2}, the analysis of Section \ref{sec-emhd} suggests that it is sufficient to show
\begin{equation}\label{psiphi}
\begin{split}
&\|\Psi(V, H)\|_{\mathcal X}+\|\Phi(V, H)\|_{\mathcal Y}\lesssim \lambda^2+ \|V\|_{\mathcal X}^2+ \|H\|_{\mathcal Y}^2 \ \ \forall \ \ (V, H)\in (\mathcal X, \mathcal Y),\\
&\|\Psi(V_1, H_1)- \Psi(V_2, H_2)\|_{\mathcal X}+\|\Phi(V_1, H_1)-\Phi(V_2, H_2)\|_{\mathcal Y}\\
\lesssim&\ \left(\lambda + \|V_1\|_{\mathcal X}+ \|V_2\|_{\mathcal X}+ \|H_1\|_{\mathcal Y}+ \|H_2\|_{\mathcal Y}\right) \left(\|V_1-V_2\|_{\mathcal X}+\|H_1-H_2\|_{\mathcal Y}\right)\\
& \ \ \ \forall \ \ (V_1, H_1)\in (\mathcal X, \mathcal Y), (V_2, H_2)\in (\mathcal X, \mathcal Y).
\end{split}
\end{equation}
Note that the estimates regarding the highest derivative term $\Phi^3$ are obtained in Section \ref{sec-emhd} and the estimates of the pure fluid term 
$\Psi^1$ are established in Section \ref{sec-nse}. The term $\Psi^2$ has a similar quadratic structure as $\Phi^3$ and it is in lower order; hence the estimates of $\Psi^2$ to establish (\ref{psiphi}) are guaranteed. Among $\Phi^1$ and $\Phi^2$ we only need to estimate one of them, say $\Phi^1$. That is the task in the following.

We start with the $\mathcal H^{\frac{7}{2}-2\alpha}(\mathbb T^3)$ norm of $V\otimes H$ from $\Phi^1$,
\begin{equation}\notag
\begin{split}
&\left\| (1-\Delta)^{\frac74-\alpha}\int_0^t e^{-(t-\tau)(-\Delta)^\alpha} \nabla\cdot (V(x, \tau)\otimes H(x,\tau))\, d\tau \right\|_{L^2_x}\\
%=&\ \left\| \int_0^t e^{-(t-\tau)(-\Delta)^\alpha}(-\Delta)^{\frac12} (1-\Delta)^{\frac74-\alpha}\frac{\nabla\cdot}{(-\Delta)^{\frac12}} (V(x, \tau)\otimes H(x,\tau))\, d\tau \right\|_{L^2_x}\\
\lesssim&\  \int_0^t \left\|e^{-(t-\tau)(-\Delta)^\alpha}(-\Delta)^{\frac12} (1-\Delta)^{\frac74-\alpha} (V(x, \tau)\otimes H(x,\tau)) \right\|_{L^2_x}\, d\tau \\
\lesssim&\  \int_0^t \left\|e^{-(t-\tau)|\xi|^{2\alpha}}|\xi|^{\frac92-2\alpha} \mathcal F\left(V(x, \tau)\otimes H(x,\tau)\right) \right\|_{L^2_{\xi}}\, d\tau.
\end{split}
\end{equation}
%where we used Plancherel's theorem in the last step. 
In view of H\"older's inequality and Lemma \ref{le-heat1} we infer
\begin{equation}\notag
\begin{split}
& \left\|e^{-(t-\tau)|\xi|^{2\alpha}}|\xi|^{\frac92-2\alpha} \mathcal F\left(V(x, \tau)\otimes H(x,\tau)\right) \right\|_{L^2_{\xi}}\\
\lesssim&\ \left\|e^{-(t-\tau)|\xi|^{2\alpha}}|\xi|^{\frac92-2\alpha}\right\|_{L^{\frac{2pq}{2p+2q-pq}}_{\xi}} \left\|\mathcal F\left( V(x, \tau)\otimes H(x,\tau)\right) \right\|_{L^{\frac{pq}{pq-p-q}}_{\xi}}\\
\lesssim&\ (t-\tau)^{-\frac{9}{4\alpha}+1-\frac{3}{2\alpha}\left(\frac{p+q}{pq}-\frac12\right)} \left\| V(x, \tau)\otimes H(x,\tau) \right\|_{L^{\frac{pq}{p+q}}_{x}}\\
\lesssim&\ (t-\tau)^{-\frac{3}{2\alpha}+1-\frac{3}{2p\alpha}-\frac{3}{2q\alpha}}  \left\|V \right\|_{L^{q}_{x}} \left\| H \right\|_{L^{p}_{x}}
\end{split}
\end{equation}
followed by applying H\"older's inequality and Lemma \ref{le-beta}
\begin{equation}\notag
\begin{split}
&\int_0^t (t-\tau)^{-\frac{3}{2\alpha}+1-\frac{3}{2p\alpha}-\frac{3}{2q\alpha}}  \left\|V \right\|_{L^{q}_{x}} \left\| H \right\|_{L^{p}_{x}}\, d\tau\\
= &\ \int_0^t (t-\tau)^{-\frac{3}{2\alpha}+1-\frac{3}{2p\alpha}-\frac{3}{2q\alpha}} \tau^{-\beta-\gamma} (\tau^\gamma\left\|V \right\|_{L^{q}_{x}})(\tau^\beta \left\| H \right\|_{L^{p}_{x}})\, d\tau\\
\lesssim&\ \left(\int_0^t (t-\tau)^{\left(-\frac{3}{2\alpha}+1-\frac{3}{2p\alpha}-\frac{3}{2q\alpha}\right)\frac{1}{1-\beta-\gamma}} \tau^{-(\beta+\gamma)\frac{1}{1-\beta-\gamma}}\, d\tau\right)^{1-\beta-\gamma}\\
&\cdot \|V\|_{L^{(\gamma,\frac{1}{\gamma})}_tL^q_x} \|H\|_{L^{(\beta,\frac{1}{\beta})}_tL^p_x}.
\end{split}
\end{equation}
The application of Lemma \ref{le-beta} is justified by the assumptions on the parameters (\ref{para}) which imply
\[\left(\frac{3}{2\alpha}-1+\frac{3}{2p\alpha}+\frac{3}{2q\alpha}\right)\frac{1}{1-\beta-\gamma}+\frac{\beta+\gamma}{1-\beta-\gamma}=1,\]
and 
\[0<\left(\frac{3}{2\alpha}-1+\frac{3}{2p\alpha}+\frac{3}{2q\alpha}\right)\frac{1}{1-\beta-\gamma}<1, \ \ 0<\frac{\beta+\gamma}{1-\beta-\gamma}<1.\]
Hence we have
\begin{equation}\label{est-psiphi1}
\begin{split}
&\left\|\int_0^t e^{-(t-\tau)(-\Delta)^\alpha} \nabla\cdot (V(x, \tau)\otimes H(x,\tau))\, d\tau\right\|_{\mathcal H^{\frac{7}{2}-2\alpha}(\mathbb T^3)}\\
\lesssim&\ \left\|V \right\|_{L^{(\gamma, \frac{1}{\gamma})}_t L^q_{x}}\left\|H \right\|_{L^{ (\beta, \frac{1}{\beta})}_t L^p_{x}}.
\end{split}
\end{equation}
Continuing to the $L^{(\beta, \frac{1}{\beta})}_t L^p_{x}$ norm we have
\begin{equation}\notag
\begin{split}
&\left\|t^\beta \int_0^t e^{-(t-\tau)(-\Delta)^\alpha}\nabla\cdot (V(x, \tau)\otimes H(x,\tau))\, d\tau \right\|_{L^{\frac{1}{\beta}}_t L^p_{x}}\\
\lesssim &\ \left\|t^\beta \int_0^t \left\|e^{-(t-\tau)(-\Delta)^\alpha} (-\Delta)^{\frac12} (V(x, \tau)\otimes H(x,\tau))\right\|_{L^p_{x}}\, d\tau \right\|_{L^{\frac{1}{\beta}}_t}
\end{split}
\end{equation}
and %Applying H\"older's inequality and Lemma \ref{le-heat1} yields
\begin{equation}\notag
\begin{split}
&\left\|e^{-(t-\tau)(-\Delta)^\alpha} (-\Delta)^{\frac12} (V(x, \tau)\otimes H(x,\tau))\right\|_{L^p_{x}}\\
\lesssim&\ \left\|e^{-(t-\tau)|\xi|^{2\alpha}} |\xi| \mathcal F\left((V(x, \tau)\otimes H(x,\tau))\right)\right\|_{L^{\frac{p}{p-1}}_{\xi}}\\
\lesssim&\ \left\|e^{-(t-\tau)|\xi|^{2\alpha}} |\xi|\right\|_{L^{q}_{\xi}}\left\| \mathcal F\left((V(x, \tau)\otimes H(x,\tau))\right)\right\|_{L^{\frac{pq}{pq-p-q}}_{\xi}}\\
\lesssim&\ (t-\tau)^{-\frac1{2\alpha}-\frac{3}{2q\alpha}}\left\|V(x, \tau)\otimes H(x,\tau)\right\|_{L^{\frac{pq}{p+q}}_{x}}\\
\lesssim&\ (t-\tau)^{-\frac1{2\alpha}-\frac{3}{2q\alpha}}\left\|V\right\|_{L^{q}_{x}}\left\|H\right\|_{L^{p}_{x}}.
\end{split}
\end{equation}
It follows that
\begin{equation}\notag
\begin{split}
&\left\| \int_0^t e^{-(t-\tau)(-\Delta)^\alpha}\nabla\cdot (V(x, \tau)\otimes H(x,\tau))\, d\tau\right\|_{L^{(\beta, \frac{1}{\beta})}_t L^p_{x}}\\
\lesssim &\ \left\|\int_0^t t^\beta(t-\tau)^{-\frac1{2\alpha}-\frac{3}{2q\alpha}}\|V\|_{L^{q}_{x}}\|H\|_{L^{p}_{x}} \, d\tau \right\|_{L^{\frac{1}{\beta}}_t}\\
\lesssim &\ \left\| \int_0^{\frac{t}{2}} t^\beta (t-\tau)^{-\frac1{2\alpha}-\frac{3}{2q\alpha}}\|V\|_{L^{q}_{x}}\|H\|_{L^{p}_{x}} \, d\tau \right\|_{L^{\frac{1}{\beta}}_t}\\
&+ \left\|\int_{\frac{t}{2}}^t t^\beta(t-\tau)^{-\frac1{2\alpha}-\frac{3}{2q\alpha}}\|V\|_{L^{q}_{x}}\|H\|_{L^{p}_{x}} \, d\tau \right\|_{L^{\frac{1}{\beta}}_t}.
\end{split}
\end{equation}
 By Young's inequality we have
\begin{equation}\notag
\begin{split}
& \left\| \int_0^{\frac{t}{2}} t^\beta (t-\tau)^{-\frac1{2\alpha}-\frac{3}{2q\alpha}}\|V\|_{L^{q}_{x}}\|H\|_{L^{p}_{x}} \, d\tau \right\|_{L^{\frac{1}{\beta}}_t}\\
\lesssim&\ \left\| \int_0^{\frac{t}{2}}  (t-\tau)^{\beta-\frac1{2\alpha}-\frac{3}{2q\alpha}}\|V\|_{L^{q}_{x}}\|H\|_{L^{p}_{x}} \, d\tau \right\|_{L^{\frac{1}{\beta}}_t}\\
\lesssim&\ \left(\int_0^{\frac{t}{2}} (t-\tau)^{(\beta-\frac{1}{2\alpha}-\frac{3}{2q\alpha})(\frac{1}{1-\beta-2\gamma})}\, d\tau\right)^{1-\beta-2\gamma}\\
&\cdot\left(\int_0^{\frac{t}{2}}\left\|H\right\|^{\frac{1}{2\beta}}_{L^{p}_{x}}\, d\tau \right)^{2\beta}\left(\int_0^{\frac{t}{2}}\left\|V\right\|^{\frac{1}{2\gamma}}_{L^{q}_{x}}\, d\tau \right)^{2\gamma}\\
\lesssim&\ \|H\|_{L^{\frac{1}{2\beta}}_t L^{p}_{x}} \|V\|_{L^{\frac{1}{2\gamma}}_t L^{q}_{x}}
\end{split}
\end{equation}
where we used the fact
\[(\beta-\frac{1}{2\alpha}-\frac{3}{2q\alpha})(\frac{1}{1-\beta-2\gamma})=-1\]
thanks to (\ref{para}), and hence $\int_0^{\frac{t}{2}} (t-\tau)^{-1}\, d\tau \lesssim 1$.

Resorting to H\"older's inequality and Hardy-Littlewood-Sobolev lemma for the parameters satisfying (\ref{para}) we deduce
\begin{equation}\notag
\begin{split}
& \left\| \int_{\frac{t}{2}}^t t^\beta (t-\tau)^{-\frac1{2\alpha}-\frac{3}{2q\alpha}}\|V\|_{L^{q}_{x}}\|H\|_{L^{p}_{x}} \, d\tau \right\|_{L^{\frac{1}{\beta}}_t}\\
\lesssim&\ \left\| \int_0^{\frac{t}{2}}  (t-\tau)^{-\frac1{2\alpha}-\frac{3}{2q\alpha}}\tau^\beta \|V\|_{L^{q}_{x}} \|H\|_{L^{p}_{x}}\, d\tau \right\|_{L^{\frac{1}{\beta}}_t}\\
\lesssim&\ \left\| \left(\int_{\frac{t}{2}}^t\tau \|H\|_{L^{p}_{x}}^{\frac{1}{\beta}} \, d\tau\right)^{\beta}\left(\int_{\frac{t}{2}}^t  (t-\tau)^{(-\frac1{2\alpha}-\frac{3}{2q\alpha})\frac{1}{1-\beta}} \|V\|_{L^{q}_{x}}^{\frac{1}{1-\beta}}\, d\tau\right)^{1-\beta} \right\|_{L^{\frac{1}{\beta}}_t}\\
\lesssim&\ \|H\|_{L^{(\beta, \frac{1}{\beta})}_t L^{p}_{x}}\left\| \int_{\frac{t}{2}}^t  (t-\tau)^{(-\frac1{2\alpha}-\frac{3}{2q\alpha})\frac{1}{1-\beta}} \|V\|_{L^{q}_{x}}^{\frac{1}{1-\beta}}\, d\tau \right\|^{1-\beta}_{L^{\frac{1-\beta}{\beta}}_t}\\
\lesssim&\  \|H\|_{L^{(\beta, \frac{1}{\beta})}_t L^{p}_{x}} 
\left\| \|V\|_{L^{q}_{x}}^{\frac{1}{1-\beta}} \right\|_{L^{\frac{1-\beta}{2\gamma}}_t}^{1-\beta}\\
%\lesssim&\  \|H\|_{L^{(\beta, \frac{1}{\beta})}_t L^{p}_{x}} \left(\int_{\frac{t}{2}}^t\left\|H\right\|^{(\frac{1}{1-\beta})(\frac{1-\beta}{2\beta})}_{L^{p}_{x}}\, d\tau \right)^{2\beta}\\
\lesssim&\ \|H\|_{L^{(\beta, \frac{1}{\beta})}_t L^{p}_{x}}\|V\|_{L^{\frac{1}{2\gamma}}_t L^{q}_{x}}.
\end{split}
\end{equation}
We conclude from the estimates above that
\begin{equation}\label{est-psiphi2}
\begin{split}
&\left\| \int_0^t e^{-(t-\tau)(-\Delta)^\alpha}\nabla\cdot (V(x, \tau)\otimes H(x,\tau))\, d\tau \right\|_{L^{(\beta, \frac{1}{\beta})}_t L^p_{x}}\\
\lesssim&\ \left\|H \right\|_{L^{\frac{1}{2\beta}}_tL^p_x}\left\|V \right\|_{L^{\frac{1}{2\gamma}}_tL^q_x}+\left\|H \right\|_{L^{(\beta, \frac{1}{\beta})}_t L^p_{x}}\left\|V \right\|_{L^{\frac{1}{2\gamma}}_tL^q_x}.
\end{split}
\end{equation}

In the end, we estimate the norm of $L^{\frac{1}{2\beta}}_t L^p_{x}$ as
\begin{equation}\notag
\begin{split}
&\left\| \int_0^t e^{-(t-\tau)(-\Delta)^\alpha} \nabla\cdot (V(x, \tau)\otimes H(x,\tau))\, d\tau \right\|_{L^{\frac{1}{2\beta}}_t L^p_{x}}\\
\lesssim &\ \left\| \int_0^t \left\|e^{-(t-\tau)(-\Delta)^\alpha} (-\Delta)^{\frac12} (V(x, \tau)\otimes H(x,\tau))\right\|_{L^p_{x}}\, d\tau \right\|_{L^{\frac{1}{2\beta}}_t}
\end{split}
\end{equation}
and %Similarly as before, we deduce by applying H\"older's inequality and Lemma \ref{le-heat1} 
\begin{equation}\notag
\begin{split}
& \left\|e^{-(t-\tau)(-\Delta)^\alpha} (-\Delta)^{\frac12} (V(x, \tau)\otimes H(x,\tau))\right\|_{L^p_{x}}\\
\lesssim &\ \left\|e^{-(t-\tau)|\xi|^{2\alpha}} |\xi| \mathcal F\left(V(x, \tau)\otimes H(x,\tau)\right)\right\|_{L^{\frac{p}{p-1}}_{\xi}}\\
\lesssim &\ \left\|e^{-(t-\tau)|\xi|^{2\alpha}} |\xi|\right\|_{L^{q}_{\xi}} \left\|\mathcal F\left(V(x, \tau)\otimes H(x,\tau)\right)\right\|_{L^{\frac{pq}{pq-p-q}}_{\xi}}\\
\lesssim &\ (t-\tau)^{-\frac{1}{2\alpha}-\frac{3}{2q\alpha}} \left\|V(x, \tau)\otimes H(x,\tau)\right\|_{L^{\frac{pq}{p+q}}_{x}}\\
\lesssim &\ (t-\tau)^{-\frac{1}{2\alpha}-\frac{3}{2q\alpha}} \|H\|_{L^{p}_{x}}\|V\|_{L^{q}_{x}}.
\end{split}
\end{equation}
Again, the last two inequalities together with Hardy-Littlewood-Sobolev lemma imply
\begin{equation}\label{est-psiphi3}
\begin{split}
&\left\| \int_0^t e^{-(t-\tau)(-\Delta)^\alpha} \nabla\cdot (V(x, \tau)\otimes H(x,\tau))\, d\tau \right\|_{L^{\frac{1}{2\beta}}_t L^p_{x}}\\
\lesssim&\ \left\| \int_0^t (t-\tau)^{-\frac{1}{2\alpha}-\frac{3}{2q\alpha}} \|H\|_{L^{p}_{x}} \|V\|_{L^{q}_{x}}\, d\tau \right\|_{L^{\frac{1}{2\beta}}_t}\\
\lesssim&\ \left\|  \|H\|_{L^{p}_{x}} \|V\|_{L^{q}_{x}} \right\|_{L^{\frac{1}{2\beta+2\gamma}}_t}\\
\lesssim&\  \left\|H \right\|_{L^{\frac{1}{2\beta}}_tL^p_x} \left\|V \right\|_{L^{\frac{1}{2\gamma}}_tL^q_x}.
\end{split}
\end{equation}
Immediately from (\ref{est-psiphi1})-(\ref{est-psiphi3}) we obtain
\begin{equation}\notag
\begin{split}
&\left\| \int_0^t e^{-(t-\tau)(-\Delta)^\alpha} \nabla\cdot (V(x, \tau)\otimes H(x,\tau))\, d\tau \right\|_{\mathcal Y}\\
\lesssim&\  \left\|V \right\|_{\mathcal X}^2+ \left\|H \right\|_{\mathcal Y}^2.
\end{split}
\end{equation}
The estimates for $V\otimes B_{f^\omega}$, $u_{g^\omega}\otimes H$ and $u_{g^\omega}\otimes B_{f^\omega}$ from $\Phi^1$ can be achieved in an analogy way. Thus we claim the first inequality of (\ref{psiphi}) is justified. As explained in Section \ref{sec-emhd}, the second inequality of (\ref{psiphi}) can be established similarly as the first one.

\bigskip

\section*{Acknowledgement}
The author is partially supported by the NSF grants DMS-1815069 and DMS-2009422, and the von Neumann Fellowship at the IAS.

\bigskip

%{\color{blue} Submission to journal: Annals of PDE, Nonlinearity, ARMA, Communications in PDE, }

%\Endrefs

\begin{thebibliography}{XX}


%\bibliographystyle{plain}
%\bibliography{NS-stability}
%\references {999}

\bibitem{ADFL}
M. Acheritogaray, P. Degond, A. Frouvelle and J-G. Liu.
\newblock {\em Kinetic formulation and global existence for the Hall-Magnetohydrodynamic system}.
\newblock Kinetic and Related Models, 4: 901--918, 2011.


%\bibitem{Bha}
%A. Bhattacharjee.
%\newblock {\em Impulsive magnetic reconnection in the Earth's magnetotail and the solar corona}.
%\newblock Ann. Rev. Astron. Astrophys., Vol. 42: 365--384, 2004.

%\bibitem{Bif}
%L. Biferale.
%\newblock {\em Shell models of energy cascade in turbulence}.
%\newblock Annu. Rev. Fluid Mech., 35: 441468, 2003.

%\bibitem{Bis94}
%D. Biskamp.
%\newblock {\em Cascade models for magnetohydrodynamic turbulence}.
%\newblock Phys. Rev. E50: 2702--2711, 1994.

%\bibitem{Bis1}
%D. Biskamp.
%\newblock {\em Magnetic reconnection in plasmas}.
%\newblock Cambridge University Press, 2000.

%\bibitem{Bis2}
%D. Biskamp.
%\newblock {\em Magnetohydrodynamic turbulence}.
%\newblock Cambridge University Press, 2003.

\bibitem{Bou96}
J. Bourgain.
\newblock {\em Invariant measures for the 2D defocusing nonlinear Schr\"odinger equation}.
\newblock Comm. Math. Phys., 176: 421-445, 1996.

%\bibitem{BDLIS}
%T. Buckmaster, C. De Lellis, P. Isett, and L. Sz\'ekelyhidi.
%\newblock {\em Anomalous dissipation for $1/5$-H\"older Euler flows}.
%\newblock Ann. of Math., Vol. 182, No. 1: 127-172, 2015.

%\bibitem{BDLS}
%T. Buckmaster, C. De Lellis, and L. Sz\'ekelyhidi.
%\newblock {\em Dissipative Euler flows with Onsager-critical spatial regularity}.
%\newblock Comm. Pure Appl. Math., Vol. 69 No. 9, 16131670, 2016.

%\bibitem{BDLSV}
%T. Buckmaster, C. De Lellis, L. Sz\'ekelyhidi, and V. Vicol.
%\newblock {\em Onsager's conjecture for admissible weak solutions}.
%\newblock Comm. Pure Appl. Math., https://doi.org/10.1002/cpa.21781. 2018.

%\bibitem{BV19}
%T. Buckmaster, and V. Vicol.
%\newblock {\em Convex integration and phenomenologies in turbulence}.
%\newblock EMS Surveys in Mathematical Sciences, 6(1):173--263, 2020.

%\bibitem{BV21}
%T. Buckmaster, and V. Vicol.
%\newblock {\em Convex integration constructions in hydrodynamics}.
%\newblock Bulletin of the AMS, 58(1): 1--44, 2021.

%\bibitem{BV}
%T. Buckmaster, and V. Vicol.
%\newblock {\em Nonuniqueness of weak solutions to the Navier-Stokes equation}.
%\newblock Ann. of Math., 189(1):101--144, 2019.


%\bibitem{BCV}
%T. Buckmaster, M. Colombo, and V. Vicol.
%\newblock {\em Wild solutions of the Navier-Stokes equations whose singular sets in time have Hausdorff dimension strictly less than 1}.
%\newblock arXiv: 1809.00600, 2018.

\bibitem{BT1}
N. Burq and N. Tzvetkov.
\newblock {\em Random data Cauchy theory for super-critical wave equation I: Local theory}.
\newblock Invent. Math., 173(3):449-475, 2008.

\bibitem{BT2}
N. Burq and N. Tzvetkov.
\newblock {\em Random data Cauchy theory for super-critical wave equation II: A global existence result}.
\newblock Invent. Math., 173(3):477-496, 2008.

\bibitem{CDL}
D. Chae, P. Degond and J-G. Liu.
\newblock {\em Well-posedness for Hall-magnetohydrodynamics}.
\newblock Ann. Inst. H. Poincar\'e Anal. Non Lineaire, Vol. 31: 555--565, 2014.

%\bibitem{CL}
%D. Chae and J. Lee.
%\newblock {\em On the blow-up criterion and small data global existence for the Hall-magneto-hydrodynamics}.
%\newblock J. Differential Equations, 256: 3835--3858, 2014.

%\bibitem{CS}
%D. Chae,  and M. Schonbek.
%\newblock {\em On the temporal decay for the Hall-magnetohydrodynamic equations}.
%\newblock J. Differential Equations,  Vol. 255: 3971--3982, 2013.

\bibitem{CWW}
D. Chae,  R. Wan and J. Wu.
\newblock {\em Local well-posedness for the Hall--MHD equations with fractional magnetic diffusion}.
\newblock J. Math. Fluid Mech., 17: 627-638, 2015.

%\bibitem{CWeng}
%D. Chae and S. Weng.
%\newblock {\em Singularity formation for the incompressible Hall-MHD equations without resistivity}.
%\newblock Ann. I. H. Poincar\'e-AN, Vol. 33: 1009--1022, 2016.

%\bibitem{CW}
%D. Chae and J. Wolf.
%\newblock {\em On partial regularity for the 3D non-stationary Hall magnetohydrodynamics equations on the plane}.
%\newblock Comm. Math. Phys., Vol. 354: 213--230, 2017.

%\bibitem{CD-Kol}
%A. Cheskidov and M. Dai.
%\newblock {\em Kolmogorov's dissipation number and the number of degrees of freedom for the 3D Navier-Stokes equations}.
%\newblock Proceedings of the Royal Society of Edinburg, Section A, Vol. 149, Issue 2: 429--446, 2019.

%\bibitem{CD}
%A. Cheskidov and M. Dai.
%\newblock {\em Regularity criteria for the 3D Navier-Stokes and MHD equations}.
%\newblock arXiv:1507.06611, 2015.

%\bibitem{CDK}
%A. Cheskidov, M. Dai and L. Kavlie.
%\newblock {\em Determining modes for the 3D Navier-Stokes equations}.
%\newblock Physica D: Nonlinear Phenomena, Vol. 374-375: 1--9, 2018.

%\bibitem{CF}
%A. Cheskidov and S. Friedlander.
%\newblock {\em The vanishing viscosity limit for a dyadic model}.
%\newblock Physica D, 238:783--787, 2009.

%\bibitem{CFP1}
%A. Cheskidov, S. Friedlander, and N. Pavlovi\'c.
%\newblock {\em Inviscid dyadic model of turbulence: the fixed point and Onsager's conjecture}.
%\newblock J. Math. Phys., 48 (6): 065503, 16, 2007.

%\bibitem{CFP2}
%A. Cheskidov, S. Friedlander, and N. Pavlovi\'c.
%\newblock {\em An inviscid dyadic model of turbulence: the global attractor}.
%\newblock Discrete Contin. Dyn. Syst., 26 (3): 781--794, 2010.

%\bibitem{CL}
%A. Cheskidov and X. Luo.
%\newblock {\em Sharp nonuniqueness for the Navier-Stokes equations}.
%\newblock arXiv:2009.06596, 2020.


%\bibitem{CSreg}
%A. Cheskidov and R. Shvydkoy.
%\newblock {\em A unified approach to regularity problems for the 3D Navier-Stokes and Euler equations: the use of Kolmogorov's dissipation range}.
%\newblock J. Math. Fluid Mech., Vol.16, Issue 2: 263--273, 2014.

%\bibitem{CheS}
%A. Cheskidov and R. Shvydkoy.
%\newblock {\em Euler equations and turbulence: analytical approach to intermittency}.
%\newblock SIAM J. Math. Anal., 46 (1): 353--374, 2014.


%\bibitem{CS}
%A. Cheskidov and R. Shvydkoy.
%\newblock {\em Ill-posedness of basic equations of fluid dynamics in Besov spaces}.
%\newblock Mathematics Subject Classification,  2000.

\bibitem{CO}
J. Colliander and T. Oh.
\newblock {\em Almost sure well-posedness of the cubic nonlinear Schr\"odinger equation below $L^2(\mathbb T)$}.
\newblock Duke Math. J., 161(3): 367-414, 2012.

%\bibitem{CDD}
%M. Colombo, C. De Lellis, and L. De Rosa.
%\newblock {\em Ill-posedness of Leray solutions for the hypodissipative Navier-Stokes equations}.
%\newblock Comm. Math. Phys., Vol.362, No.2, 659688, 2018.

%\bibitem{CLT}
%P. Constantin, B. Levant, and E.Titi.
%\newblock {\em Analytic study of the shell model of turbulence}.
%\newblock Physica D: Nonlinear Phenomena, 219 (2): 120--141, 2006.

%\bibitem{CF}
%P. Constantin and C. Foias.
%\newblock {\em Naiver-Stokes Equations}.
%\newblock Chicago Lectures in Mathematics, The University of Chicago Press, 1988.

%\bibitem{CC}
%A. C\'ordoba and D. C\'ordoba.
%\newblock {\em A maximum principle applied to quasi-geostrophic equations}.
%\newblock Comm. Math. Phys., 249(3): 511--528, 2004.

%\bibitem{Dai-20}
%M. Dai.
%\newblock {\em Blow-up of a dyadic model with intermittency dependence for the Hall MHD}.
%\newblock  Physica D: Nonlinear Phenomena, Vol. 428: 133066, 2021.

%\bibitem{Dai1}
%M. Dai.
%\newblock {\em Local well-posedness of the Hall-MHD system in $H^s(\mathbb R^n)$ with $s>\frac n2$}.
%\newblock  Mathemaische Nachrichten. DOI: 10.1002/mana.201800107, 2019.

\bibitem{Dai2022}
M. Dai.
\newblock {\em Almost sure existence of global weak solutions for supercritical electron MHD}.
\newblock arXiv:2201.08161, 2022.

\bibitem{Dai2}
M. Dai.
\newblock {\em Local well-posedness for the Hall-MHD system in optimal Sobolev spaces}.
\newblock Journal of Differential Equations, 289: 159-181, 2021.

%\bibitem{Dai3}
%M. Dai.
%\newblock {\em Propagation of regularity for the MHD system in optimal Sobolev space}.
%\newblock arXiv: 1707.07754, 2017.

%\bibitem{Dai18}
%M. Dai.
%\newblock {\em Non-unique weak solutions in Leray-Hopf class of the 3D Hall-MHD system}.
%\newblock arXiv: 1812.11311, 2018.

%\bibitem{Dai4}
%M. Dai.
%\newblock {\em Regularity criterion for the 3D Hall-magneto-hydrodynamics}.
%\newblock Journal of Differential Equations. Vol. 261: 573--591, 2016.

%\bibitem{D}
%M. Dai.
%\newblock {\em Regularity criterion and energy conservation for the supercritical Quasi-Geostrophic equation}.
%\newblock Journal of Mathematical Fluid Mechanics. To appear. ArXiv:1505.02293, 2015.

%\bibitem{DF}
%M. Dai and S. Friedlander.
%\newblock {\em Dyadic models for ideal MHD}.
%\newblock arXiv:2104.09440. To appear in the Journal of Mathematical Fluid Mechanics.

%\bibitem{DL}
%M. Dai and H. Liu.
%\newblock {\em Long time behavior of solutions to the 3D Hall-magneto-hydrodynamics system with one diffusion}.
%\newblock Journal of Differential Equations, Vol. 266: 7658--7677, 2019.

%\bibitem{DL-well}
%M. Dai and H. Liu.
%\newblock {\em On well-posedness of generalized Hall-magneto-hydrodynamics}.
%\newblock arXiv: 1906.02284, 2019.

\bibitem{DT}
R. Danchin and J. Tan.
\newblock {\em On the well-posedness of the Hall-magnetohydrodynamics system in critical spaces}.
\newblock Communications in Partial Differential Equations, 46(1):31-65, 2021.


%\bibitem{DSz}
%S. Daneri, and L. Sz\'ekelyhidi.
%\newblock {\em Non-uniqueness and h-principle for H\"older-continuous weak solutions of the Euler equations}.
%\newblock Arch. Ration. Mech. Anal., Vol. 224 No.2: 471--514, 2017.

%\bibitem{DLS1}
%C. De Lellis, and L. Sz\'ekelyhidi, Jr.
%\newblock {\em Dissipative continuous Euler flows}.
%\newblock Invent. Math., Vol.193 No. 2: 377--407, 2013.

%\bibitem{DLS2}
%C. De Lellis, and L. Sz\'ekelyhidi, Jr.
%\newblock {\em Dissipative Euler flows and Onsager's conjecture}.
%\newblock Journal of the European Mathematical Society,  16(7): 1467--1505, 2014.

%\bibitem{DLS3}
%C. De Lellis, and L. Sz\'ekelyhidi, Jr.
%\newblock {\em The Euler equations as a differential inclusion}.
%\newblock Ann. of Math.,  Vol.170 No.3: 1417--1436, 2009.

%\bibitem{DR}
%L. De Rosa.
%\newblock {\em Infinitely many Leray–Hopf solutions for the fractional Navier-Stokes equations}.
%\newblock Communications in Partial Differential Equations, 44(4): 335--365, 2019.

\bibitem{CD}
C. Deng and S. Cui.
\newblock {\em Random-data Cauchy problem for the Navier-Stokes equations on $\mathbb T^3$}.
\newblock J. Differential Equations, 251: 902-917, 2011.

\bibitem{Deng}
Y. Deng.
\newblock {\em Two-dimensional nonlinear Schr\"ondinger equation with random radial data}.
\newblock Anal. PDE, 5: 913-960, 2012.

%\bibitem{DN}
%V. N.  Desnyanskiy and E. A. Novikov.
%\newblock {\em Evolution of turbulence spectra toward a similarity regime}.
%\newblock Izv. Akad. Nauk SSSR, Fiz. Atmos. Okeana, 10: 127--136, 1974.


%\bibitem{DS}
%E. I. Dinaburg and Y. G. Sinai.
%\newblock {\em A quasi-linear approximation of three-dimensional Navier-Stokes system}.
%\newblock Moscow Math. J., 1: 381--388, 2001.

\bibitem{DLM}
B. Dodson, J. L\"uhrmann and D. Mendelson.
\newblock {\em Almost sure local well-posedness and scattering for the 4D cubic nonlinear Schr\"odinger equation}.
\newblock Adv. Math., 347: 619-676, 2019.

%\bibitem{DG}
%C. Doering and J. D. Gibbon.
%\newblock {\em Applied Analysis of the Navier-Stokes Equations}.
%\newblock Cambridge Texts in Applied Mathematics, Cambridge University Press, 2004.

%\bibitem{DZ}
%L. Du and T. Zhang.
%\newblock {\em Almost sure existence of global weak solutions for incompressible MHD equations in nagative-order Sobolev space}.
%\newblock J. Differential Equations, 263: 1611-1642, 2017.

%\bibitem{DS}
%E. Dumas and F. Sueur.
%\newblock {\em On the weak solutions to the Maxwell-Landau-Lifshitz equations and to the Hall-magnetohydrodynamic equations}.
%\newblock Comm. Math. Phys., 330: 1179--1225, 2014.

%\bibitem{FLS}
%D. Faraco, S. Lindberg, and L. Sz\'ekelyhidi.
%\newblock {\em Bounded solutions of ideal MHD with compact support in space-time}.
%\newblock arXiv: 1909.08678, 2019.

%\bibitem{FMRR}
%C. L. Fefferman, D. S. McCormick, J. C. Robinson and J. L. Rodrigo.
%\newblock {\em Higher order commutator estimates and local existence for the non-resistive MHD equations and related models}.
%\newblock Journal of Functional Analysis, Vol. 267: 1035--1056, 2014.

%\bibitem{FHN}
%J. Fan, S. Huang and G. Nakamura.
%\newblock {\em Well-posedness for the axisymmetric incompressible viscous Hall-magnetohydrodynamic equations}.
%\newblock Appl. Math. Lett., 26: 963--967, 2013.

%\bibitem{Fi}
%N. Filonov.
%\newblock {\em Uniqueness of the Leray-Hopf solution for a dyadic model}.
%\newblock Transactions of the American Mathematical Society, Vol. 369 (12): 8663--8684, 2017.

%\bibitem{FK}
%N. Filonov and P. Khodunov.
%\newblock {\em Non-uniqueness of Leray-Hopf solutions for a dyadic model}.
%\newblock St. Petersburg Math. J., Vol. 32: 371--387, 2021.

%\bibitem{FP}
%S. Friedlander and N. Pavlovi\'c.
%\newblock {\em Blowup in a three-dimensional vector model for the Euler equations}.
%\newblock Comm. Pure Appl. Math., 57 (6): 705--725, 2004.

%\bibitem{Fri}
%U. Frisch.
%\newblock {\em Turbulence: The Legacy of A. N. Kolmogrov}.
%\newblock Cambridge University Press, Cambridge, 1995.

%\bibitem{Galtier}
%S. Galtier.
%\newblock {\em Introduction to Modern Magnetohydrodynamics}.
%\newblock Cambridge University Press, London, 2016.

%\bibitem{Gle}
%E. B. Gledzer.
%\newblock {\em System of hydrodynamic type admitting two quadratic integrals of motion}.
%\newblock Soviet Phys. Dokl., 18: 216-217, 1973.

%\bibitem{Gr}
%L. Grafakos.
%\newblock {\em Modern Fourier analysis}.
%\newblock Second edition. Graduate Texts in Mathematics, 250. Springer, New York, 2009.

%\bibitem{Ho}
%E. Hopf.
%\newblock {\em \"Uber die Anfangswertaufgabe f\"ur die hydrodynamischen Grundgleichungen}.
%\newblock Math. Nachr., 4:213--231, 1951.

%\bibitem{Is0}
%P. Isett.
%\newblock {\em Holder continuous Euler flows with compact support in time}.
%\newblock Ph. D. Thesis, Princeton Jniversity, 2013.

%\bibitem{Is}
%P. Isett.
%\newblock {\em A Proof of Onsager's Conjecture}.
%\newblock Ann. of Math., Vol.188 No.3: 1--93, 2018.

%\bibitem{JO}
%I. Jeong and S. Oh.
%\newblock {\em On the Cauchy problem for the Hall and electron magnetohydrodynamic equations without resistivity I: illposedness near degenerate stationary solutions}.
%\newblock arXiv: 1902.02025, 2019.

%\bibitem{JS1}
%H. Jia and V. \v{S}ver\'ak.
%\newblock {\em Are the incompressible 3d Navier-Stokes equations locally ill-posed in the natural energy space?}
%\newblock J. Funct. Anal., Vol. 268(12): 3734--3766, 2015.

%\bibitem{JS2}
%H. Jia and V. \v{S}ver\'ak.
%\newblock {\em Local-in-space estimates near initial time for weak solutions of the Navier-Stokes equations and forward self-similar solutions}.
%\newblock Invent. Math., Vol. 196(1): 233--265, 2014.

%\bibitem{KP}
%N. Katz and N. Pavlovi\'c. 
%\newblock {\em Finite time blow-up for a dyadic model of the Euler equations}.
%\newblock Trans. Amer. Math. Soc., 357 (2): 695--708, 2005.

%\bibitem{KZ}
%A. Kiselev and A. Zlato\v{s}. 
%\newblock {\em On discrete models of the Euler equation}.
%\newblock Int. Math. Res. Not., 38: 2315--2339, 2005.

%\bibitem{K41}
%A. Kolmogoroff.
%\newblock The local structure of turbulence in incompressible viscous fluid for very large {R}eynold's numbers.
%\newblock {\em C. R. (Doklady) Acad. Sci. URSS (N.S.)}, 30:301--305, 1941.

%\bibitem{Lan} L. D. Landau.
%\newblock {\em On the problem of turbulence}.
%\newblock Doklady Akademii Nauk SSSR. 44: 339--342, 1944.

\bibitem{LRS}
G. Lebowitz, H. Rose, and E. Speer.
\newblock {\em Statistical mechanics of the nonlinear Schr\"odinger equation}.
\newblock J. Statist. Phys., 50(3-4): 657-687, 1988.


\bibitem{LR}
P. G. Lemari\'e-Rieusset.
\newblock {\em Recent developments in the Navier-{S}tokes problem}.
\newblock Chapman and Hall/CRC Research Notes in Mathematics, 431. Chapman  and Hall/CRC, Boca Raton, FL, 2002.

%\bibitem{Le}
%J. Leray.
%\newblock {\em Sur le mouvement d'un liquide visqueux emplissant l'espace}.
%\newblock Acta Math., 63(1):193--248, 1934.

%\bibitem{Lions}
%J.-L. Lions.
%\newblock {\em Quelques m\'ethodes de r\'esolution des probl\'emes aux limites non lin\'eaires}.
%\newblock volume 1. Dunod; Gauthier-Villars, Paris, 1969.

\bibitem{LM}
J. L\"uhrmann and D. Mendelson.
\newblock {\em Random data Cauchy theory for nonlinear wave equations of power-type on $\mathbb R^3$}.
\newblock Comm. Partial Differential Equations, 39(12): 2262-2283, 2014.

%\bibitem{LT}
%T. Luo and E. Titi.
%\newblock {\em Non-uniqueness of weak solutions to hyperviscous Navier-Stokes equations-On sharpness of J.-L. Lions exponent}.
%\newblock arXiv:1808.07595, 2018.

%\bibitem{Lu}
%X. Luo.
%\newblock {\em Stationary solutions and nonuniqueness of weak solutions for the Navier-Stokes equations in high dimensions}.
%\newblock arXiv:1807.09318, 2018.

%\bibitem{LPPPV}
%V. S. L'vov, E. Podivilov, A. Pomyalov, I. Procaccia, and D. Vandembroucq.
%\newblock {\em Improved shell model of turbulence}.
%\newblock Phys. Rev. E (3) 58: 1811--1822, 1998.

%\bibitem{MB}
%A. J. Majda and A. L. Bertozzi.
%\newblock {\em Vorticity and incompressible flow}.
%\newblock Cambridge Jniversity Press, Cambridge, JK, 2001.


%\bibitem{MH}
%H. Miura and D. Hori.
%\newblock {\em Hall effects on local structure in decaying MHD turbulence}.
%\newblock J. Plasma Fusion Res., 8: 73--76, 2009.

%\bibitem{MS}
%S. Modena, and L. Sz\'ekelyhidi.
%\newblock {\em Non-uniqueness for the transport equation with Sobolev vector fields}.
%\newblock arXiv:1712.03867v3, 2017.

\bibitem{NPS}
A. R. Nahmod, N. Pavlovi\'c, and G. Staffilani.
\newblock {\em Almost sure existence of global weak solutions for super-critical Navier-Stokes equations}.
\newblock SIAM J. Math. Anal., 45(6): 3431-3452, 2013.

%\bibitem{Ob}
%A. M. Obukhov. 
%\newblock {\em Some general properties of equations describing the dynamics of the atmosphere}.
%\newblock Izv. Akad. Nauk SSSR Ser. Fiz. Atmosfer. i Okeana, 7:695--704, 1971.


%\bibitem{Og}
%T. Ogawa.
%\newblock {\em Sharp Sobolev inequality of logarithmic type and the limiting regularity condition to the harmonic heat flow}.
%\newblock SIAM J. Math. Anal., 34: 1318--1330, 2003.

%\bibitem{OY}
%K. Ohkitani and M. Yamada.
%\newblock {\em Temporal intermittency in the energy cascade process and local Lyapunov analysis in fully-developed model of turbulence}.
%\newblock Progr. Theoret. Phys., 81: 329--341, 1989.

%\bibitem{On}
%L. Onsager.
%\newblock {\em Statistical hydrodynamics}.
%\newblock Nuovo Cimento (9), 6(Supplemento, 2(Convegno Internazionale di Meccanica Statistica)):279--287, 1949.


%\bibitem{PM}
%M. Polygiannakis and X. Mossas.
%\newblock {\em A review of magneto-vorticity induction in Hall- MHD plasmas}.
%\newblock Plasma Phys. Control \& Fusion, 43: 195--221, 2001.

\bibitem{Po}
O. Pocovnicu.
\newblock {\em Almost sure global well-posedness for the energy-critical defocusing nonlinear wave equation on $\mathbb R^d$, $d=4$ and $5$}.
\newblock J. Eur. Math. Soc., 19(8): 2521-2575, 2017.

%\bibitem{Sch}
%V. Scheffer.
%\newblock {\em An inviscid flow with compact support in space-time}.
%\newblock J. Geom. Anal., 3(4):343--401, 1993.

%\bibitem{SJ}
%D. Shalybkov and V. Jrpin.
%\newblock {\em The Hall effect and the decay of magnetic fields}.
%\newblock Astron. Astrophys., 685--690, 1997.

%\bibitem{Sh1}
%A. Shnirelman.
%\newblock {\em On the nonuniqueness of weak solution of the Euler equations}.
%\newblock Comm. Pure Appl. Math., 50(12):1261--1286, 1997.

%\bibitem{Sh2}
%A. Shnirelman.
%\newblock {\em Weak solutions with decreasing energy of incompressible Euler equations}.
%\newblock Comm. Math. Phys., 210(3):541--603, 2000.

%\bibitem{Tao}
%T. Tao.
%\newblock {\em A quantitative formulation of the global regularity problem for the periodic Navier-Stokes equation}.
%\newblock Dyn. Partial Differ. Eq., 4(4): 293-302, 2007.

%\bibitem{Tao}
%T. Tao.
%\newblock {\em Finite time blowup for an averaged three-dimensional Navier-Stokes equation}.
%\newblock J. Amer. Math. Soc., Vol. 29: 601--674, 2016.

%\bibitem{Tem}
%R. Temam.
%\newblock {\em Navier-Stokes Equations: Theory and Numerical Analysis}.
%\newblock AMS Chelsea, Providence, Rhode Island, 2000.

%\bibitem{Wal}
%F. Waleffe.
%\newblock {\em On some dyadic models of the Euler equations}.
%\newblock Proc. Amer. Math. Soc., 134 (10): 2913--2922, 2006.

%\bibitem{Wa}
%M. Wardle.
%\newblock {\em Star formation and the Hall effect}.
%\newblock Astrophys. Space Sci., 292: 317--323, 2004.


\end{thebibliography}
\end{document}